\documentclass[12pt]{article}
\pdfoutput=1

\usepackage{arxiv}
\usepackage[utf8]{inputenc} 
\usepackage[T1]{fontenc}    
\usepackage{hyperref}       
\usepackage{url} 
\usepackage{booktabs}       

\usepackage{nicefrac}       
\usepackage{microtype}      
\usepackage{lipsum}
\usepackage{graphicx}
\usepackage{float}

\usepackage{amsmath}
\usepackage{amsthm}
\usepackage{amsfonts}       
\usepackage{amssymb,eurosym,latexsym}

\usepackage{fullpage}
\usepackage{subfigure}
\usepackage{mathtools}
\usepackage{dsfont}
\usepackage{tikz}
\usepackage{booktabs}
\usepackage{multirow}
\usepackage{multicol}
\usepackage{bigints}
\usepackage{appendix}
\usepackage{xcolor}

\restylefloat{figure}
\restylefloat{table}

\usetikzlibrary{calc}
\usetikzlibrary{arrows}
\usetikzlibrary{shapes.geometric}

\newtheorem{theorem}{Theorem}
\newtheorem{corollary}[theorem]{Corollary}
\newtheorem{example}{Example}
\newtheorem{lem}[theorem]{Lemma}

\newtheorem{proposition}{Proposition}
\newtheorem{definition}{Definition}

\theoremstyle{remark}
\newtheorem*{remark}{Remark}

\newcommand{\Max}{\mathrm{max}}
\newcommand{\Min}{\mathrm{min}}
\newcommand{\X}{{\cal{X}}}
\newcommand{\A}{{{\cal{A}}}}

\newcommand{\St}{{\cal{S}}}
\newcommand{\Esp}{\mathbb{E}}

\allowdisplaybreaks

\usepackage[linesnumbered,ruled,vlined]{algorithm2e}

\SetCommentSty{mycommfont}
\SetKwInput{KwInput}{Input}                
\SetKwInput{KwOutput}{Output}              

\usetikzlibrary{shapes,decorations,arrows,calc,arrows.meta,fit,positioning}
\tikzset{
    -Latex,auto,node distance =1 cm and 1 cm,semithick,
    state/.style ={ellipse, draw, minimum width = 0.7 cm},
    point/.style = {circle, draw, inner sep=0.04cm,fill,node contents={}},
    bidirected/.style={Latex-Latex,dashed},
    el/.style = {inner sep=2pt, align=left, sloped}
}

\title{Optimal control policies for resource allocation in the Cloud: comparison between Markov decision process and heuristic approaches}

\author{
Thomas Tournaire 
\thanks{Nokia Bell Labs France, \texttt{thomas.tournaire@nokia.com}.}\,
\thanks{Samovar, Telecom SudParis.}
\and
Hind Castel-Taleb 
\footnotemark[2]
\and
Emmanuel Hyon
\thanks{
Sorbonne Universit{\'e}, UPMC Univ Paris 06, 
CNRS, LIP6 UMR 7606.}\,
\thanks{Université Paris Nanterre.}
}

\begin{document}
\maketitle

\begin{abstract}
We consider an auto-scaling technique in a cloud system where virtual machines hosted on a physical node are turned on and off depending on the queue's occupation (or thresholds), in order to minimise a global cost integrating both  energy consumption and performance. 
We propose several efficient optimisation methods to find threshold values minimising this global cost:  local search heuristics coupled with aggregation of Markov chain and with queues approximation techniques to reduce the execution time and improve the accuracy.  
The second approach tackles the problem with a Markov Decision Process (MDP) for which we proceed to a theoretical study and provide theoretical comparison with the first approach. We also develop structured MDP algorithms integrating hysteresis properties.  We show that MDP algorithms (value iteration, policy iteration) 
and especially structured MDP algorithms outperform the devised heuristics, 
in terms of time execution and accuracy. 
Finally, we propose a cost model for a real scenario of a cloud system  to apply our optimisation algorithms and show their relevance. 
\end{abstract}

\section{Introduction}
\label{sec:introduction}

Rapid growth of the demand for computational power by scientific, business and web-applications has led to the creation of large-scale data centers which are often oversized to ensure the quality of service of hosted applications. This leads to an under-utilisation of the servers which implies important electric losses and over-consumption. Nowadays, Cloud computing requires more electric power than those of whole countries such as India or Germany \cite{report:oie,art:mastelic2015}. 
To improve the utilisation  rate of servers, some data center owners have deployed some methods implementing dynamicity of resources according to system load. Those mechanisms, called ``{\it autoscaling}'' \cite{art:lorido2014} are based on activation and deactivation of Virtual Machines (VMs) \cite{art:benoit} according to the workload. 

Nevertheless, it is crucial to analyse both the energy consumption and the performance of the system to find the policy that adapts resource allocation to the  demand.  Unfortunately, these two measures are inversely proportional, which motivates researchers to evaluate them simultaneously via a unique global cost function. 
In order to adapt resources, some queueing systems modulate the service capacity according to queue occupancy, these are hysteresis models 
\cite{art:asghari14}.
They allow not to activate and deactivate the servers too frequently when the load is varying.
This makes it possible to correctly adapt variable resources according to demand  by means of thresholds \cite{art:serfozo,art:Lui99}.
Therefore, multi-server queueing systems working with  thresholds-based policies and verifying hysteresis properties have been suggested to efficiently manage the number of active VMs  \cite{conf:SPS15,art:asghari14,art:Tou19}.
The other advantage that makes hysteresis policies so appealing 
is their ease of implementation. This is why, they are a key component of
the auto-scaling systems for the widespread cloud architectures Kubernetes
for docker components, and Azure or Amazon Web Services \cite{rap:aws} for
virtual machines.
Henceforth, there is a great interest of studying the computation of hysteresis policy since threshold values can be plugged into autoscaling systems that are implemented in the major cloud architectures.

In this work, we consider a cost-aware approach. We propose a model with a
global expected cost considering different costs: associated with the
performance requirements defined in a Service Level Agreement (SLA) and
associated with energy consumption. 
Two modelling approaches are addressed: static optimisation problem or
dynamic control.
The static optimisation problem expresses the problem
by using the average value of the costs computed by
stationary distribution of a CTMC (Continuous Time Markov Chains).
Such an approach is often used in the field of
inventory management for base stock policies \cite{art:warsing2013}.
On the other hand, dynamic control expresses the problem under a MDP (Markov Decision Process) model and computes the optimal policy.
 
We propose new scalable and efficient algorithms 
for minimisation of the global cost.
For static optimisation the work \cite{art:Tou19} presents and compares
several heuristics but is restricted to the static model. In the present
paper, we first improve the previous heuristics of \cite{art:Tou19} by modifying several stages and present new algorithms for the static model.
We design completely new algorithms based on optimal dynamic control integrating  hysteresis  policies and we show that they are much faster while ensuring the optimality of the solution. 
Although they are the most widely used for threshold computation \cite{liv:song2013}, the efficiency of these two approaches are not assessed very well in the literature for the computation of single thresholds by level. Especially, this had been never done for hysteresis policies. Determining which is the best promising approach is thus a major point which is addressed here by performing numerous numerical experiments.
Furthermore, numerical results show that we can generate the thresholds corresponding to the optimized global cost in just few seconds even for large systems. This could have a significant impact for a cloud operator who wants to dynamically allocate its resources to control its financial costs.

The key contributions of this paper are as follows:
\begin{enumerate}
\item We improve the existing heuristics of \cite{art:Tou19} that makes a large review of the heuristics of the literature which solve static optimisation problems. 
We assess numerically the gain given  by these improvements;
\item   We provide a theoretical study of the multichain properties of the Markov decision process. This shows that some usual algorithms solving MDP models have no convergence guarantee.  We still provide an algorithm that solves dynamic control model by considering hysteresis assumption in the MDP;
\item We made a theoretical analysis between the two approaches and give some insights which explain why the static optimisation problem is suboptimal.
Numerical studies show that the dynamic control approach strongly outperforms the other one for optimality and running time criteria;
\item We develop and analyse a financial cost model. It takes into account prices of VM instantiations from cloud providers as well as energy consumption of VMs.  A presentation of minimised costs
for a problem based on this concrete cloud model is done.
\end{enumerate} 

The remainder of this paper is organised as follows. Section \ref{sec:RelWork} briefly reviews the related works. In Section \ref{sec:cloudModel}, we provide a treatment to unify the different definitions of hysteresis policies coming from
different models. Then, we describe the cloud system, with the
queueing model. We detail the cost function used to express the expected costs in terms of performance and energy consumption for such models. 
In Section \ref{sec:MCmodel}, we present the Markov chain approach
with  the decomposition and aggregation method, and
the optimisation algorithms  based on local search heuristics. We also present a Markov decision process model and adapted 
algorithms in Section
\ref{sec:MDPmodel}. Section \ref{sec:concreteModel} presents 
a concrete cost model. 
Numerical experiments  for the  comparison of the
algorithms  are discussed in Section \ref{sec:numericalExp}. 
In the conclusion, we finally discuss about achieved results. Some comments about further researches are given.


\section{Literature review}
\label{sec:RelWork}
\subsection{Energy and Performance Management in the Cloud}
In recent years, the increase in the energy consumption has remained a crucial problem for large-scale computing. 
In a data center, the server power consumption can be divided into static and dynamic parts. The static part (which does not vary with
workload) represents the energy consumed by a server 
when it is idle,
while the dynamic cost  depends on the current usage.
In \cite{art:orgerie}, they 
define a power-aware model  to estimate the dynamic 
part of energy cost 
for a VM  of a given size, this model keeps the philosophy of 
the pay as you go model but it is based on energy consumption.  

As the static part represents a high part of the overall
energy consumed by the server nodes, therefore, shutting unused
physical resources that are idle  lead to non-negligible energy savings.
Two main approaches of physical server resource management have been proposed to improve the energy efficiency: shutdown or
switching on servers or VMs which is referred as dynamic power management \cite{art:benoit},  
and  scaling of the CPU performance which is
referred as Dynamic Voltage and Frequency Scaling \cite{art:Krzy18}.
Shutdown strategies (considered here) are often combined with 
consolidation algorithms that gather the load on 
few servers to favour the shutdown of the others. 
So, managing energy by switching on or switching off virtual machines is an intuitive and widespread manner to save energy.
Yet, as quoted in \cite{art:benoit}, coarse techniques 
of shutdown are, most often, not the appropriate solution to
achieve energy reduction. 
Indeed, shutdown policies suffer from energy and time losses when
switching off and switching on take longer than the idle period.

\subsection{Control Management for Queueing Models}
In order to represent the problems with activations and deactivations of
virtual machines, server farm models have been 
proposed \cite{art:adan2014,Gandhi2010,art:mitrani2013managing}. 
Usually, these server farm models are modeled with multi-server
queueing systems \cite{Artalejo2005,art:ardagna2014}.
Although server farms with multi-server queues allow a fairly fine
representation of the dynamicity induced by virtualization, these models
do not address issues related neither to the internal network
nor to the VM placement in it.  
All VMs are instead considered as parallel resources. This makes these 
models suitable for studying simple nodes of several servers in the 
cloud.
If queuing models allow us to easily compute performance metrics, the
decision making for switching on or switching off the VM  requires an additional step which remains a key point. The computation of the optimal actions has led to a large field of researches and methods. 
Dynamic control and especially Markov decision Processes appear 
to be the main direct method. 

\subsection{Markov Decision Process and Hysteresis Policies}
Since the seminal work of Mc Gill in 1969 (with an optimal control model) 
and that of Lippman (with a Markov decision process model) 
numerous works have been devoted to similar multi-server queue models 
using Markov decision processes (see \cite{art:teghem85} and references therein for the oldest, and more recently \cite{proc:YCNH11} and
\cite{art:Lee2014} to quote 
just a few). Unfortunately, not all of them received rigorous treatment and the study of unichain or multichain property is often ignored. 
It appeared very early, in the work of Bell for a $M/M/2$ model (see 
\cite{art:teghem85}), that the optimal policy in such models has a special form and is called \emph{hysteresis}.
In hysteresis policies, servers are powered up when the number of jobs in the system is sufficiently high and are powered down when that number is
sufficiently low. 
More precisely, 
activations  and deactivations of servers are ruled by 
sequences of \emph{different} forward and reverse thresholds. 
The concept of hysteresis can also be applied to the control of service rates.

Researches on hysteresis policies are twofold: 
first exploration of the conditions that insure the
optimality of  hysteresis policy. 
Hence, Szarkowicz et al. \cite{art:szarkowicz} showed the optimality of
hysteresis policies in a $M/M/S$ queueing model with a 
control of the vacations of the servers. For models with control of the service
rates, the proofs are made in Hipp \cite{art:hipp84} or Serfozo \cite{art:serfozo}. The second axis studies the computation of the 
threshold values.

\subsection{Threshold calculation methods}
For hysteresis models, the calculation of optimal thresholds
received less attention in the past. Also here two major trends appeared.
The computation of the optimal policy can be done by means of
adapted usual dynamic programming algorithms in which the
structured policies properties are plugged. 
A similar treatment has been addressed for routing in 
\cite{art:NinoMora2019}. The alternate way is similar to the single threshold research for base-stock policies which is very common
in inventory management. A local search \cite{art:randa2019,art:kranenburg} is used to explore the optimal thresholds. The computation of expected measures associated with a set of thresholds is complex \cite{art:braglia2019}. It requires 
the computation of the stationary distribution either by standard numerical methods see e.g. \cite{art:Wu2014} or after a Markovian analysis with simpler and faster computations (e.g. \cite{art:warsing2013} uses iterative methods to compute the stationary distribution).
In  \cite{liv:song2013}, Song claims that, for single threshold
models, the second approach dealing with stationary distribution
computations are
generally more effective than the MDP approach.
Such a comparison has not been performed yet for hysteresis
especially since few works implement a MDP algorithm with a
structured policy.

As cloud systems are modelled by multi-dimensional systems,  defined on very large state spaces,  then
the stationary distribution computation is difficult.  
Fortunately, the computation of the performance measures of
hysteresis multi-server systems has been already studied in the literature. 
Different efficient resolution methods have been developed.
Among the most significant works, 
we quote the work of Ibe and Keilson
\cite{art:ibekeilson95} refined in Lui and Golubchik \cite{art:Lui99}.
Both solve the model by partitioning the state space in disjoint sets to aggregate the Markov chain. 
Exhaustive comparisons of the resolution methods are made in \cite{conf:epewKandi17}:
closed-form solution of the stationary distribution, 
partition of the state space, and matrix geometric methods applied on  QBD (Quasi birth and death) processes are studied.
It was noticed that partitioning is the most suited.
Furthermore, for optimisation, the objective cost function 
being non-linear and non-convex makes this problem very complex and
there is currently no exact method to solve this problem. 
The work  \cite{art:kranenburg} presents three heuristics for searching single thresholds in inventory models that have been
adapted for hysteresis in \cite{art:Tou19}. These approximate
heuristics require the computation of the invariant distribution
and the cost for numerous threshold values which requires very high computation times. 
For a server farm model with activation of a single reserve
threshold, Mitrani \cite{art:mitrani2013managing}
uses fluid approximation to compute the activation thresholds
and \cite{art:Wu2014} uses genetic algorithm for optimisation and matrix geometric method for the stationary distribution computation.

\section{Cloud Model}
\label{sec:cloudModel}
In this section, we present the cloud system,
denoted \textit{Cloud Model}, and modelled  by a controlled multi-server queueing model, where  arrivals of requests are processed by servers. 
The servers in the queueing system represent logical units in the cloud. 
Since we consider IaaS clouds (Infrastructure as a Service) which provide
computing resources as a service, the logical units are either virtual
machines or docker components.
We propose a mathematical analysis of the queueing model in order to derive performance as well as energy consumption measures.  

\subsection{Controlled multi-server queue}
\label{sec:controlledQueue}
We have following assumptions for the model:
\begin{enumerate}
\item Arrivals of requests follow a Poisson process of rate
$\lambda$, and service times of all VMs are independent of arrivals
and independent of each other.
Moreover, they are i.i.d. and follow an exponential distribution
with an identical rate $\mu$;
\item Requests enter in the queue, which has a finite capacity  $B$, 
and are served by one of the $K$ servers. 
The service discipline is supposed  FIFO (First In First Out).
\end{enumerate} 
A customer is treated by a server
as soon as this server becomes idle and the server is active.
Servers can be turned on and off by the controller. 
When the server is turned on it 
becomes active while it becomes inactive when it is turned off. 

We define $\mathcal{S}$ the state space, where  $\mathcal{S}=\{0,1,\ldots,B\} \times \{1,\ldots,K\}$. Any state $x \in \mathcal{S}$ is such that $x=(m,k)$ where  
$m$ represents the number of requests in the system, 
and $k$ is the number of operational servers (or active servers, this number can also be seen as a service level). 
We define $\mathcal{A}=\{0,\ldots,K\}$, be the set of actions, where action $a$ denotes the number of servers to be activated. 
With this system, are associated two kinds of costs that a cloud provider encounters: 
\begin{enumerate}
\item Costs corresponding  to the  performance of the system,  for the control of the  service  quality defined  in the  SLA: as costs ($C_H$) per unit of time for holding requests in the queue or
instantaneous costs ($C_R$) for losses of requests.   
\item Costs corresponding to  the use of resources (operational and 
energy consumption): as costs for using a VM per unit of time  ($C_S$) 
and instantaneous costs for activating/deactivating ($C_A$ and $C_D)$.
\end{enumerate} 
We define for the system a global  cost containing  both performance and energy consumption costs. 
We have the following objective function we want to minimise:
\vspace{-0.5em}
\begin{equation}\label{eq:globalobjective}
\bar{C} =\lim_{T \rightarrow \infty}\frac{1}{T} \, \Esp \Biggl\{ \int_0^T C_t\left(X(t),A(t)\right)dt \Biggr\}
\end{equation}
where $A(t)$ is the action taken at time $t$ (it is possible that nothing was performed) and $C_t(X(t),A(t))$ is the cost that is charged over time when the state is $X(t)$ and action
$A(t)$ is performed.
This problem can be solved  using 
dynamic programming algorithms or by computing the corresponding
thresholds values for turning on and off the VMs.

\subsection{Models of hysteresis policies}
\label{sec:hysteresisPolicy}
A decision rule is a mapping from some
information set to some action. A policy is 
a sequence of decision rules $\eta=(q_0,q_1,q_2,\ldots)$.
The most general set of policies is that of history-dependent
randomised policies, but the classical results on average
infinite-horizon, time-homogeneous Markovian optimal control
\cite{liv:put1994} allow us to focus on stationary Markov Deterministic Policies. 
Such policies are characterised by a single, deterministic decision
rule which maps the current state to an action. 
We thus consider the  mapping $q$ from $\mathcal{S} $ to $\mathcal{A}$
such that $q(x)=a$. 

All along the paper we restrict our attention with a special 
form of policies: the \emph{hysteresis} policies. 
Nevertheless, there is relatively few homogeneities 
between definitions of hysteresis  policies in
the literature. Indeed, it can refer to policies defined with
double thresholds (especially when 
$K=2$), or to a restrictive definition when Markov chains are used.
We follow the works of \cite{art:hipp84,art:serfozo} to present an
unified treatment of hysteresis.

\paragraph{Hysteresis policies with multiple activations}
We assume in this part that the decision rule is a mapping 
from the state to a number of active servers 
$q(m,k)=k_1$ with $k_1 \in [1,\ldots,K]$.
Several servers can be activated or deactivated 
to pass from $k$ to $k_1$ active servers.

\begin{definition}[Double threshold policies]
\label{def:doubleLevel}
We call \emph{double threshold} policy a stationary  
policy such that 
the decision rule $q(m,k)$ is increasing in both of its arguments
and is defined by a sequence of thresholds such that 
for any $k$ and $k_1$ in $[1,K]$ we have:
\[ q(m,k) = k_1 \text{ when } \quad \ell_{k_1}(k) \leq m < \ell_{k_1+1}(k) \, ,\]
where $\ell_{k_1}(k)= \min \{m :  q(m,k) \geq  k_1\}$.
This minimum is $\infty$ if the set is empty. 
For all $k$, we also fix $\ell_{K+1}(k)=\infty$ and $\ell_{1}(k)=0$ (since
at least one server must be active).
\end{definition}

A monotone hysteresis policy is a special case of double threshold policy.
\begin{definition}[Monotone Hysteresis polices \cite{art:hipp84}]
\label{def:monotoneHysteresis}
A policy is a monotone hysteresis policy  
if it is a double threshold policy and moreover if there 
exist two sequences of integers $l_k$ and $L_k$ such that
\begin{gather*}
l_k = \ell_{k}(K) = \ell_{k}(K-1) = \ldots = \ell_{k}(k) \\
L_k= \ell_{k}(1) = \ell_{k}(2) = \ldots = \ell_{k}(k-1) \, ,
\end{gather*}
with $l_k \leq L_k     \  \text{ for } \  k = 1, \ldots, K+1$  ;  $l_k \leq L_{k+1} \  \text{ for } \ k=1, \ldots, K \,$ ; $l_1 = L_1 = 0$, $l_{K+1} = L_{K+1} = \infty$. And if, for all  $(m,k) \in \mathcal{S}$,  
\begin{equation*}
 q(m,k) =
\begin{cases}
q(m,k-1)  & \text{ if } m < l_k  \text{ and }  k > 1 \\
k         & \text{ if } l_k \leq m < L_{k+1} \\
q(m,k+1) & \text{ if } m \geq L_{k+1} \text{ and } k < m 
\end{cases} \, .
\end{equation*}
\end{definition}
The thresholds $l_{k}$ can be seen as the queue levels at which some
servers should be deactivated and the $L_{k+1}$ are the analogous
activation points.
Roughly speaking, the difference between a double threshold and a monotone
hysteresis policy lies in the fact that some thresholds toward a level 
are identical in hysteresis.

\begin{definition}[Isotone Hysteresis polices \cite{art:hipp84}]]
\label{def:isotoneHysteresis}
A policy is an isotone policy if it is a monotone hysteresis policy and if $0= l_1 \leq l_2 \leq \ldots \leq l_{K+1} = \infty \text{ and } 
0 = L_1 \leq L_2 \leq \ldots \leq L_{K+1} = \infty$.
\end{definition}

\begin{example}
In Figure \ref{fig:iso}, we represent an isotone policy. The number
indicates the number of servers that should be activated in each state. 
The bold line means that the number of activated servers is the same than 
the decision and then no activation or deactivation have to be performed.
\begin{figure}[hbtp]
\begin{center}
    \includegraphics[scale=0.4]{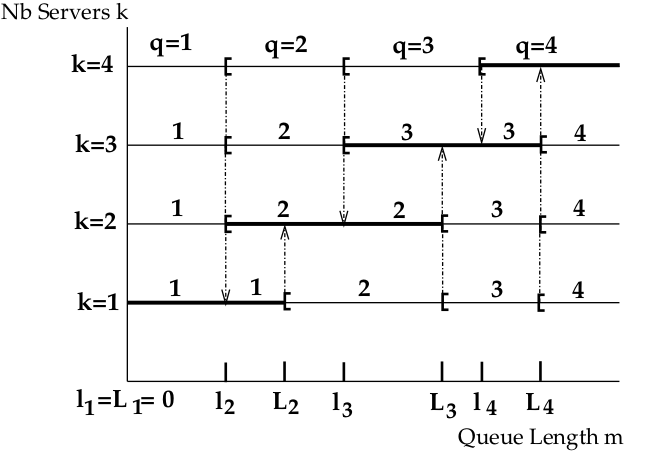}
\end{center}
\caption{An isotone hysteresis policy.}\label{fig:iso}
\end{figure}
\end{example}

\begin{proposition}[Optimality of monotone hysteresis policies \cite{art:szarkowicz}]
In the multi-server queueing model presented in Section \ref{sec:cloudModel} and for which there is activation/deactivation costs, working costs and holding costs, Szarkowicz and al. showed that monotone hysteresis policies are optimal policies. 
\end{proposition}

However, the presence of a rejection cost, as we assume here, is not
considered in the assumptions of \cite{art:szarkowicz} and there is no
proof about the optimality of hysteresis policies for the model studied
here.

\begin{remark}
There is an alternate way to express the policy by giving the number of
servers to activate or deactivate instead of giving directly the number of servers that should be active. 
\end{remark}

\paragraph{Hysteresis policies with single VM activation}
We focus now on models in which we can only operate one machine at a time.
Now the decision rule indicates an activation or deactivation.
Therefore, the action space is now $\mathcal{A} = \{ -1, 0, 1 \}$. 
For a state $(m,k)$, when 
the action $q \in \mathcal{A}$ is $-1$ then we deactivate a server,
when the action is $0$ then the number of active servers remains
identical, when the action is $1$ then we activate a server.

It could be noticed that, for this kind of model, double threshold policies and hysteresis policies coincide. 
Indeed, the decision now relates to activation 
(resp. deactivation) and no
longer to the number of servers to activate (resp. deactivate). 
There exist only two thresholds by level $k$:
$L_{k+1}$ to go to level $k+1$, and $l_k$ to go to level $k-1$.
There are no other levels that can be reached from $k$. 
For example, in Figure \ref{fig:iso} all decisions smaller 
than the level
are replaced by $-1$ while all decisions larger than the level are
replaced by $1$.
Definition \ref{def:isotoneHysteresis} remains unchanged, but Definition \ref{def:monotoneHysteresis} should be rephrased in:

\begin{definition}[Monotone hysteresis policy]\label{def:hysteresisSVM}
A policy is a monotone hysteresis policy if it is a stationary policy
such that the decision rule $q(m,k)$ is increasing in $m$ and decreasing
in $k$ and if is
defined by two sequences of thresholds $l_k$ and $L_k$ such that
for all $(m,k) \in \mathcal{S}$: 
\begin{equation*}
q(m,k) =
\begin{cases}
-1  & \text{ if } m < l_k  \text{ and }  k > 1 \\
\ 0        & \text{ if } l_k \leq m < L_{k+1} \\
\ 1 & \text{ if } m \geq L_{k+1} \text{ and } k < m
\end{cases} \, ,
\end{equation*}
with $l_k \leq L_k$ for $k = 1, \ldots, K+1$ ; 
$l_k \leq L_{k+1}$  for $k=1, \ldots, K$  and  
$l_1 = L_1 = 0 \, ; l_{K+1} = L_{K+1} = \infty.$
\end{definition}

\paragraph{Hysteresis policies and Markov chain}
As presented in \cite{art:teghem85}, there exists a slightly different
model of multi-server queue with hysteresis policy which received a lot
of attention (see Ibe and Keilson \cite{art:ibekeilson95} or
Lui and Golubchik \cite{art:Lui99} and references therein). 
This model is still a multi-server queueing system but is no longer a
controlled model: the transitions between levels are also
governed by sequences of prefixed thresholds this is why it is called an
hysteresis model.
It is built to be easily represented by a Markov chain.
The differences between the controlled model and this Markov chain model are detailed in Section \ref{sec:modelDifferences}.

\begin{remark}
For convenience, deactivation and activation thresholds will be denoted respectively by $R$ and $F$ in the Markov chain model, while they will be denoted $l$ and $L$ in the MDP model.
\end{remark}

\begin{definition}[Hysteresis policy \cite{art:Lui99}]\label{def:hystGolub}
A $K$-server threshold-based queueing system with hysteresis 
is defined by a sequence  $F=[F_{1}, F_{2}, \ldots, F_{K-1}]$ of activation 
thresholds 
and a sequence $[R_{1}, R_{2}, \ldots, R_{K-1}]$ of deactivation thresholds.
For $1\leq k < K$, the threshold $F_k$  makes the system goes from level $k$ to 
level $k+1$ when a customer arrives 
with $k$ active servers and  $F_k$ customers in the system.
Conversely, the threshold $R_k$  makes the system goes from level $k+1$ to level $k$ when a customer leaves with $k+1$ active servers and $R_{k}+1$ customers in the system.

Furthermore, we assume that $F_1<F_2<\ldots<F_{K-1}\leq K$,  
$1 \leq R_1<R_2<\ldots<R_{K-1}$, and $R_k<F_k$, $\forall \, 1 \leq k \leq K-1$. We denote the vector that gathers the two threshold vectors $F$ and $R$ by $[F,R]$.
\end{definition}

\begin{remark}
It can be noticed that,
no server can remain idle all the time here since the threshold values are
bounded, whereas in Definition \ref{def:hysteresisSVM}  when a threshold
is infinite the server remains inactive. Hence the hysteresis policy
presented in \cite{art:ibekeilson95} can be seen as a restricted
version of isotone policies given in \ref{def:hysteresisSVM}.
Furthermore, the inequalities being strict in \cite{art:Lui99}, 
the hysteresis of
Definition \ref{def:hystGolub} is a very specific case
of hysteresis of Definition \ref{def:hysteresisSVM} that we call \emph{strictly isotone}.
\end{remark}

\vspace{-1em}
\section{Hysteresis policies and Markov chain approach}
\label{sec:MCmodel}
This section is devoted to the study of the policies defined in
Definition \ref{def:hystGolub} and their aggregation properties.
   
\subsection{Hysteresis queueing model} \label{sec:mc}
In the hysteresis model of Definition \ref{def:hystGolub}
thresholds are fixed before the system works. 
For a set  $[F,R]$ of fixed thresholds, the stochastic model is represented by a Continuous-Time Markov Chain (CTMC).
The chain is denoted by  $\{X(t)\}_{t \geq 0}$ and the 
state space is denoted by $\X \subset \St$.
Each state $(m,\,k)$ in $\X$ is 
such that $m$ is the number of requests in the 
system and $k$ is the number of active servers.  
Thus, the state space is given by:
\begin{equation}
\X= \! \{(m,\,k) ~ | ~    0 \leq m \leq F_{1},  \text{ if } k=1  \, , R_{k-1}+1 \leq m \leq F_{k}, \text{if } 1 < k < K \, , R_{K-1}+1 \leq m \leq B,  \text{if } k=K
\} .
\end{equation}
The transitions between states are described by:
\vspace{-0.4cm}
\begin{align*}
(m,k)   \rightarrow ~ &  (\min\{B,m+1\},k), \, \text{with rate } \lambda \, , 
\ \text{ if } m < F_{k} \, ; \\
     \rightarrow ~ &  (\min\{B,m+1\},\min\{K,\,k+1\}), \text{with rate } \lambda \, , \ \text{ if } m=F_k \, ; \\
    \rightarrow ~ &  (\max\{0,m  - 1\},k), \, \text{with rate} \, \mu \! \cdot \!  \min\{m,k\}, \ \text{if } m  > R_{k-1} + 1  ; \\
    \rightarrow ~&  (\max\{0,m  -1\},\max\{1,k -1\}) \, \text{with rate} \,  \mu  \min\{k,m\} \, , \
                 \text{if} \,  m  =  R_{k - 1}  + 1  \, .    
\end{align*}

\begin{figure}[hbtp]
    \subfigure[]{\includegraphics[width=0.495\textwidth]{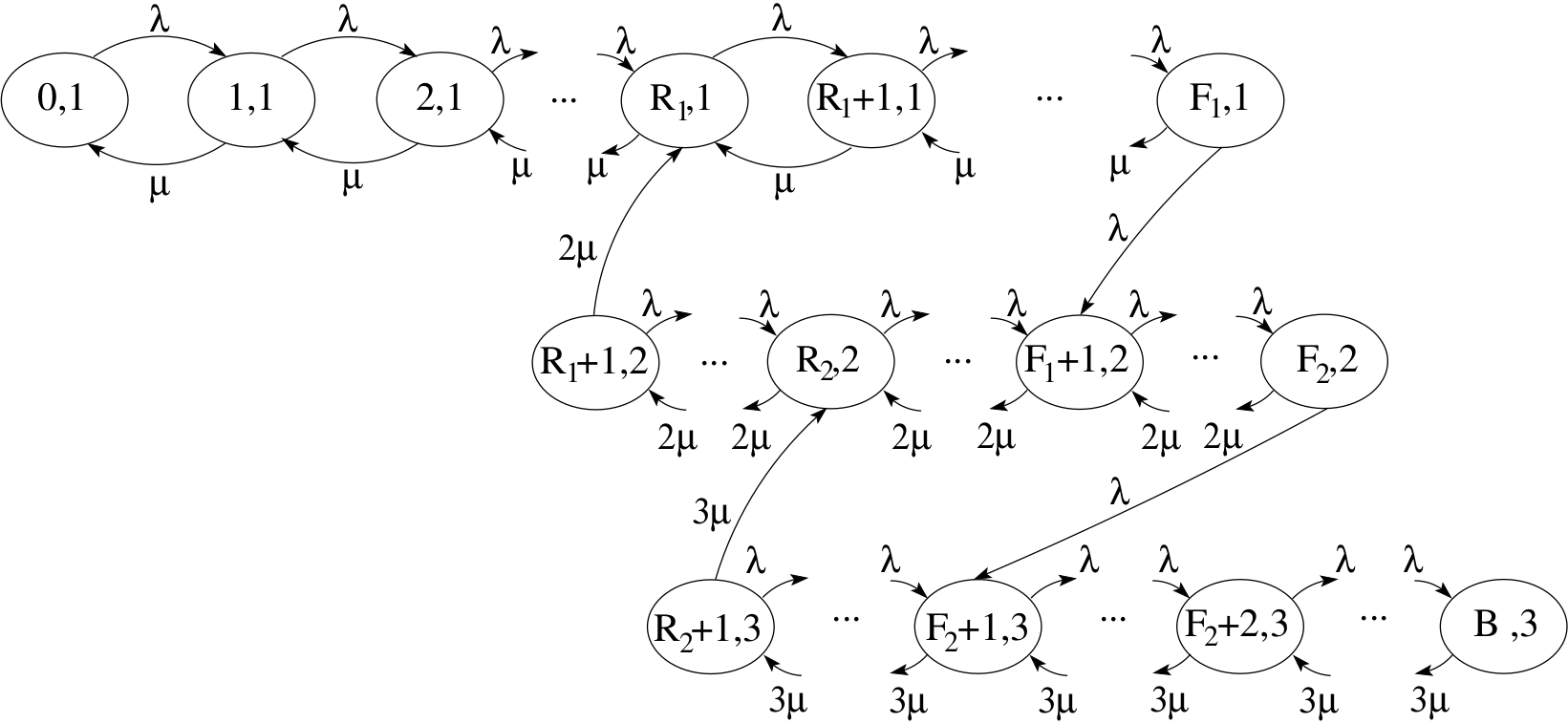}} 
    \subfigure[]{\includegraphics[width=0.495\textwidth]{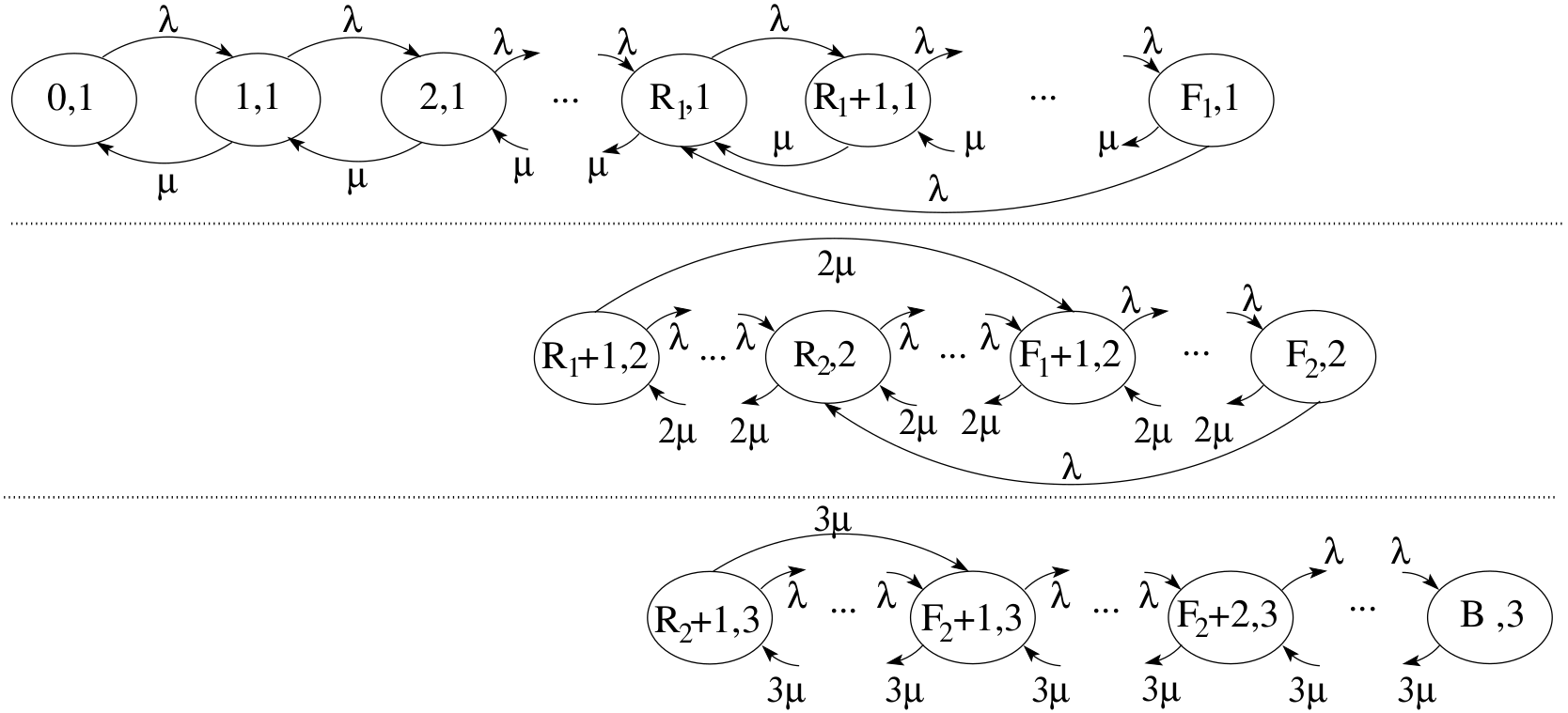}} 
    \caption{(a) Markov chain for K=3 , \hspace{20pt} (b) Micro-chains aggregation}
    \label{fig:Chaine}
\end{figure}

In  Fig. \ref{fig:Chaine}.a,  we give an example of these
transitions for a maximum number of requests in the system equal to $B$ and a number of levels $K$ equal to three.

\subsubsection{Costs}\label{sec:defmodelcost}
The global cost defined in Equation \eqref{eq:globalobjective} can be rewritten for fixed $[F,R]$. We have:
\begin{equation}
\overline{C}_{[F,R]}^{\pi} ~ = ~ \sum \limits_{k=1}^K ~~ \sum \limits_{m=R_{k-1} + 1}^{F_k} ~ \pi(m,k) \cdot \overline{C}(m,k) \, ,
\label{eqcout}
\end{equation}
where $\pi(m,k)$ represents the stationary probability in state $(m,k)$.
For a given state of the environment (number of requests, number of activated VMs), the cloud provider pays a cost:
\vspace{-1em}
\begin{multline}
\overline{C} (m , k) =C_H \cdot m +  C_S \cdot k +  C_A \cdot \lambda \cdot \mathds{1}_{m = F_k,\; k < K} + C_D \cdot \mu \cdot \min\{m,k\} \cdot \mathds{1}_{m = R_{k-1}+ 1, \; 2 \leq k \leq K } \\
+ C_R \cdot \lambda \cdot \mathds{1}_{m = B, \; k=K  } \, .
\end{multline}

\subsection{Aggregation methods to compute the global cost} 
\label{sec:aggregationMCmodel}
Solving the Markov chain quickly becomes complex when the  number of levels $K$ increases. So we propose 
to apply the SCA (Stochastic Complement Analysis) method \cite{art:Lui99}, based on a decomposition and an aggregation of the Markov chain. 
This approach allows a simplified computation of the  exact
stationary probability.  
The aim is to separate the $K$ levels of the considered Markov
chain into $K$ independent sub-chains, called \emph{micro-chains}
and then to solve each of them independently. 
We then build the aggregated 
Markov chain which has $K$ states such that each state represents an 
operating level. This chain is called \emph{macro-chain}.

\paragraph{Details about \emph{micro-chains} and \emph{macro-chain}} \label{agreg}
For any level $k$, such that $1 \leq k \leq K$,
the micro-chain is built from the initial Markov chain by keeping 
the transitions inside the same level.
As for the transitions between states of different levels, 
they are deleted and replaced by new transitions inside the
\emph{micro-chains} from a leaving state to an entering state of the initial chain. 
Let us illustrate this in Fig.~\ref{fig:Chaine}.a for initial chain
and in Fig.~\ref{fig:Chaine}.b for the micro-chains.
Inside level 2, the transitions with rate $\lambda$ and those with rate $2\mu$ are kept, while the transitions between level 2 and level 1 and  those between level 2 and level 3 are removed. 
The deleted transitions between levels 2 and 3 are replaced by a transition inside level $2$ between $(F_2,2)$ and $(R_2,2)$ with rate $\lambda$, and inside level $3$ by a transition between $(R_2+1,3)$ and $(F_2+1,3)$ with rate $3\mu$. 
Similarly for transitions between levels 2 and 1, then the process is generalised for any levels.  

The \emph{macro-chain}, is represented by a $K$ states Markov
chain, such that each state is a level $k$ 
(with $1 \leq k  \leq K$). 
The \emph{macro-chain} is a birth-death process for which computing the stationary probability is well-known.  
Let $\Pi$ be the stationary distribution of the \emph{macro-chain} 
and $\pi_k$ be the stationary distribution of 
the \emph{micro-chain} of level $k$ (with $1 \leq k \leq K$). 
We denote by $\pi$ the stationary distribution of the initial Markov 
chain. From the SCA method, the probability distribution $\pi$  is computed as follows:
\vspace{-0.5em}
\begin{equation}
\forall ~ (m,k) \in \mathbb{S} ,\; \pi(m,k) ~ = ~ \Pi(k) \cdot \pi_k(m)  \, .
\label{eqDistribSCA}
\end{equation}

Since, as pointed out in \cite{conf:epewKandi17}, 
the use of closed formulas suffers from numerical instability,
the solving of Markov chains (\emph{micro} and \emph{macro}) is done numerically using the \emph{power method} implemented in the \textit{marmoteCore} software \cite{conf:ajmMarmote}.
Each \emph{micro-chain} is solved separately, 
next the \emph{macro-chain}, and lastly, 
with Eq.~\ref{eqDistribSCA}, the initial Markov chain.  
The distributions obtained with the aggregation process are
compared with those obtained by direct computations on the whole chain. The difference is always smaller than $10^{-8}$ which is the precision we choose for our computations.

\paragraph{Aggregated Costs} 
The aggregation framework described above, allows to compute 
the global cost $\overline{C}_{[F,R]}^{\pi}$,
defined by Eq.~\eqref{eqcout}, in an aggregated form 
based on the costs per level computed with \emph{micro-chains}. 

\begin{theorem}[Aggregated Cost Function]
\label{AgregationCout}
Let $[F,R]$ be fixed thresholds, then the expected cost $\overline{C}_{[F,R]}^{\pi}$ of the system can be decomposed as:
\begin{displaymath}
\overline{C}_{[F,R]}^{\pi}  ~ = ~ \sum \limits_{k=1}^K ~ \Pi(k) \cdot \overline{C}(k) \, , \quad \text{ with } \, \, \overline{C}(k) =\sum \limits_{m=R_{k-1} + 1}^{F_k}  \pi_k(m) \cdot \overline{C}(m,k) \, ,
\end{displaymath}
where $\overline{C}(k)$ represents the expected cost of level $k$. 
\end{theorem}

\begin{proof}
We use Eq. (\ref{eqDistribSCA}) in Eq. (\ref{eqcout}) to get  $\overline{C}_{[F,R]}^{\pi} = ~ \sum \limits_{k=1}^K ~~ \sum \limits_{m=R_{k-1} + 1}^{F_k} ~ \pi_k(m) \cdot \Pi(k) \cdot \overline{C}(m,k) \, . $
Because $\Pi(k)$ does not depend on $m$, it is removed 
from the second sum. This yields the average aggregated costs per level
weighted by the \emph{macro-chain} stationary distribution.  
\begin{displaymath}
\overline{C}_{[F,R]}^{\pi} ~  = ~ \sum \limits_{k=1}^K ~ \Pi(k) ~ \cdot \sum \limits_{m=R_{k-1} + 1}^{F_k} ~ \pi_k(m) \cdot \overline{C}(m,k) ~ = ~ \sum \limits_{k=1}^K ~ \Pi(k) \cdot \overline{C}(k) \, .
\end{displaymath}
\vspace{-1em}
\end{proof}

Theorem \ref{AgregationCout} suggests an efficient approach for the global cost computation. 
Indeed, instead of having computations on a
multidimensional Markov chain, we have computations over several one-dimensional Markov chains.
We take advantage of this aggregated expression to improve the computation
of the expected cost when only a single threshold changes.

\begin{corollary}
\label{cor1}
Let $[F,R]$ be a fixed vector of thresholds. It is assumed that the \emph{micro-chains} and the costs per level associated 
with $[F,R]$ are
already computed. The modification of a threshold $F_k$ or $R_k$ in
$[F,R]$ only impacts the levels $k$ and $k+1$ 
in the Markov chain.
Therefore, the computation of the new average cost only requires a new
computation of $\pi_k, \pi_{k+1}, \overline{C}(k), \overline{C}(k+1)$ and $\Pi$.
\end{corollary}
The intuition about this corollary can be deduced
from the transition graph in Fig. \ref{fig:Chaine}.a. 
When modifying $F_2$, it only impacts the distributions
$\pi_2$ and $\pi_3$. We thus have to compute new $\overline{C}(2)$ and
$\overline{C}(3)$, while the last expected cost $\overline{C}(1)$
remains unchanged. The details of the proof are below.

\begin{proof}
Recall that we define the load of the system by 
$\rho = \lambda / \mu$.
Using the balance equations in \cite{art:Lui99}, 
we exhibit the impact
generated by the modification of activation and deactivation
thresholds on the stationary probability distributions of the
\emph{micro-chains} first
and then of the \emph{macro-chain}.

Impact on the \emph{micro-chains}: \ 
Owing to $m \in [R_{k-1} + 1 , F_k]$ for all $k \in \{1, \ldots, K \}$, 
we can therefore express the stationary probability of each state in the level $k$ 
by: $\pi_k(m) ~ = ~ \pi_k(R_{k-1}+1) \cdot \gamma_m^k $ where:
\begin{displaymath}
\pi_k(R_{k-1}+1)=  \Biggl( \sum_{m=R_{k-1}+1}^{F_k} \gamma_m^k \Biggr)^{-1} \, ,
\end{displaymath} 
\vspace{-2em}

and
\begin{equation*}
 \gamma_m^k = 
\begin{cases}
\sum \limits_{j=0}^{m - R_{k-1} - 1} \left(\dfrac{\rho}{k}\right)^j \, , 
     ~ \mbox{ if } R_{k-1} + 1 \leq m \leq R_k \, ;\\
\sum \limits_{j=0}^{m - R_{k-1} - 1} \left(\dfrac{\rho}{k}\right)^j  ~ -      
     \gamma^k_{F_k} \cdot   \sum \limits_{j=1}^{m - R_k} \left(\dfrac{\rho}{k}\right)^j \, , 
     ~ \mbox{ if }  R_k + 1 \leq m \leq F_{k-1} + 1 \, ;\\
\sum \limits_{j=m - F_{k-1} - 1}^{m - R_{k-1} - 1} \left(\dfrac{\rho}{k}\right)^j  ~  -   ~         \gamma^k_{F_k} \cdot \sum \limits_{j=1}^{m - R_k} \left(\dfrac{\rho}{k}\right)^j \, ,
     ~ \mbox{ if }  F_{k-1} + 2 \leq m \leq F_k - 1 \, ;\\
\rho \Bigg[ \left(\dfrac{\rho}{k}\right)^{F_k -      F_{k-1}} - \left(\dfrac{\rho}{k}\right)^{F_k - R_{k-1}+1}  \Bigg] /
\left(k + \rho \left(\dfrac{\rho}{k}\right)^{F_k - R_k} \right)  \, , 
     ~ \mbox{  if } m = F_k \, .
\end{cases}
\end{equation*}
Note that  $\pi_k(R_{k-1}+1)$  is determined through normalisation condition which states that the sum of probabilities of all states of the \emph{micro-chain} of level $k$ equals 1.

We notice, in the equations above, that the stationary probability 
formulas of $\pi_k(m)$ of levels $k$ 
only depend on the thresholds $R_{k-1}$, $R_k$, $F_{k-1}$ and $F_k$.
This shows that the variation of a threshold $R_k$ or $F_k$ has an 
impact only on levels $k$ and $k+1$ but not on the other levels. 
Moreover, since the \emph{micro-chain} of level $k$ starts from state ($R_{k-1} + 1 , k$) and ends in state ($F_k , k$), then thresholds 
from level $1$ to $k-2$ as well as thresholds from level $k+2$ to $K$ are not involved on the stationary distribution of the level $k$.

To resume, if we modify an activation threshold $F_k$, then it will 
modify states $(F_k , k)$ and $(F_k + \! 1 , k \! + \!1)$ and the two \emph{micro-chains} of level $k$ and $k+1$. 
If we modify, a deactivation threshold $R_{k-1}$, then it will modify states $(R_{k-1} + \! 1 , k)$ and $(R_{k-1} , k-1)$ 
and the two \emph{micro-chains} of levels $k-1$ and $k$.

Finally, the variation of a threshold from level $k$ 
modifies only the two \emph{micro-chains} of level $k$ and $k+1$.
Hence, to compute the cost with respect to 
Theorem \ref{AgregationCout} we only need to recalculate the stationary probabilities and thus the associated costs of these levels. The costs of the other levels are left unchanged. 

Impact on the \emph{macro-chain}: \ 
A change in a threshold value modifies some transitions of the \emph{macro-chain}. This can be seen, similarly as previously, using the balance equations: $\lambda_k ~ = ~ \lambda \cdot \pi_k(F_k)$, $\forall  ~ k ~ = ~ 1 ,  \ldots  , K-1$, and   $\mu_k ~ = ~ \mu \cdot \pi_k (R_{k-1} + 1)$, $  \forall ~ k ~ = ~ 2,\ldots,  K$. 
This gives the following formulas for the stationary probability distribution, for all $k \in \{2 ,  \ldots  , K\}$:
\begin{displaymath}
\begin{cases}
\lambda_k ~ = ~ \lambda \cdot \pi_k(F_k) 
 \qquad \qquad \ \ \forall  ~ k ~ = ~ 1 ,  \ldots  , K-1\\
 \mu_k ~ = ~ \mu \cdot \pi_k (R_{k-1} + 1) 
 \qquad \forall ~ k ~ = ~ 2,\ldots,  K
\end{cases} \, .
\end{displaymath}
\end{proof}

\subsection{Thresholds calculation using stationary distributions of Markov chain}
\label{sec:heuristicsMCmodel}

Our first approach relies on thresholds calculation using the stationary
distribution of the underlined Markov chain coupled with an optimisation
problem. We focus here on different static algorithms which compute the optimal threshold policy minimising the expected global cost.

\paragraph{The optimisation problem}
\label{sec:OptimisationMCmodel} 
For a set of fixed parameters: $\lambda$, $\mu$, $B$, $K$, and costs: $C_a$, $C_d$, $C_h$, $C_s$, $C_r$, 
the algorithms seeks to find the vector of thresholds $[F,R]$, 
such that the objective function is minimal.
The objective function is the expected cost 
$\overline{C}_{[F,R]}^{\pi}$ defined by Eq.\eqref{eqcout}.
For any solution $[F,R]$, the computation of the expected cost depends on the stationary distribution $\pi$ of the Markov chain $\{X(t)\}_{t \geq 0}$ induced by $[F,R]$.
The thresholds must verify some given constraints 
resulting from the principle of hysteresis. 
That is: $R_k < F_k$, $R_{k-1} < R_k$ and $F_{k-1} < F_k $. Also, it is assumed that servers do not work for free. That is: $k < F_k$. Therefore, we have to solve the constrained optimisation problem given by: 
\begin{equation}
\label{eqOptim}
\begin{aligned}
& \underset{F , R}{\text{Minimise}}
& & \hspace{1cm} \overline{C}_{[F,R]}^{\pi}  \\
& \text{u.c.}
& & R_i < F_i ~~~ \; i = 1, \ldots,K - 1. \\
&&&  F_1 < F_2 < \cdots < F_{K-1} <  B \\
&&&  0 \leq R_1 < R_2 < \cdots < R_{K-1}\\
&&& F_i ~ , ~ R_i \in \mathbb{N}
\end{aligned} \, .
\end{equation}

Since the cost function is non-convex \cite{art:Tou19} 
and threshold  values should be integer, 
we face up a combinatorial 
problem whose resolution time increases exponentially with the 
number of thresholds.

\subsubsection{Local search heuristics}\label{sec:Localsearch}
To solve this optimisation problem, several local search heuristics and one meta-heuristic have been developed in \cite{art:Tou19}.
Most of the heuristics are inspired by Kranenburg et al. \cite{art:kranenburg} 
which presents three different heuristics to resolve a lost sales
inventory problem with, unlike us, 
only a single threshold by level. 
Actually usual base-stock heuristics had to be adapted to the
hysteresis model as it is necessary to manage double 
thresholds now. 
The work \cite{art:Tou19} compares both the accuracy 
and execution time of these heuristics in numerical experiments.
We keep here the two heuristics \textbf{BPL}, \textbf{NLS} that have the best results and present some improvement ways.

\paragraph{\textit{Best Per Level} (\textbf{BPL})}
This algorithm first initialises $[F,R]$ with the lowest feasible value, i.e. $F_1 = 1, F_2 = 2, ... , F_{K-1} = K$ and $R_1 = 0, R_2 = 1, ... , R_{K-1} = K-1$.
Then it improves each threshold in the following order: it starts with the first activation threshold $F_1$ by testing all its possible values taking into account the hysteresis constraints. A new value of $F_1$ that improves the global cost, will replace the old one. Then, it will move on to $F_2$, $F_3$, ... , $F_{K-1}, R_1 , R_2, ... , R_{K-1}$. Once a loop is finished, it will restart again until the mean global cost is not improved anymore.

\paragraph{\textit{Neighborhood Local Search} (\textbf{NLS})}
This algorithm is the classical local search algorithm. We randomly
initialise the solution $[F,R]$. Then we generate the neighborhood $\mathcal{V}(x)$ of a current solution $x$. 
Each neighboring solution is the same as the current solution except for a shift $\pm 1$ for one of the thresholds. 
The algorithm explores all the neighborhood and returns the best solution among $\mathcal{V}$. Again, it loops the same process until the mean global cost is not improved anymore.

\subsubsection{Aggregation for local search algorithms}
Local search algorithms are based on an exploration step in a set of feasible solutions that must be evaluated. 
This evaluation consists of calculating 
$C^{\pi}_{([F,R])}$. From Section \ref{sec:aggregationMCmodel}, this
computation depends on the stationary distributions $\pi_k$ of the
\emph{micro-chains}, on the distribution $\Pi$ of the \emph{macro-chain},
and on the costs per level. 
Although, the neighboring solutions differ from the current solution only by a single change of one of their thresholds, after any modification we have to recalculate all these elements. 
This requires a huge amount of time which depends on the number of neighboring solutions tested at each iteration 
and  which increases when we consider large scale scenarios (with large $K$ and $B$).

We ensure the algorithms to avoid so many computations for each
neighboring solution with Corollary \ref{cor1}. Its use will widely
accelerate the algorithm without altering the accuracy.
Indeed, at each iteration, we only have to compute the stationary
distributions and the corresponding average costs per level of the two \emph{micro-chains} impacted by a change of a threshold. 
We no longer need to compute the distributions of the other \emph{micro-chains} since we can keep their average costs per level.
The macro-chain still need to be solved however.

Although decomposition and aggregation approaches reduce 
the number of computations to perform, 
some algorithmic tricks must be implemented in order
to ensure a proper and efficient running. Hence,
the two micro-chains of the current solution that are modified
when we study a neighborhood, have to be stored since they will be
reused many times. By this way, the execution time of aggregated
local search heuristics is reduced. This will be all the more effective when $K$ increases.


\subsubsection{Solution using queueing model approximations.}
We work here on a new heuristic which calculates a near-optimal solution very quickly. It can be considered as a method in its own
right but can be also used for the initialisation step in local search heuristics.
This heuristic is called \textit{MMK approximation} since it
uses ${M}/{M}/{k}/{B}$ results to compute  costs of these fairly close models.

\paragraph{Principle of the heuristic}
The main idea of the heuristic is that, at each level $k$, 
we compute an approximation of the mean global cost by using 
the stationary distribution of the Markov chain of a  $M/M/k/B$ queue
model instead of using the Markov chain of Section \ref{sec:mc}.
In such queues, the stationary distribution of the Markov chain is given
by a closed formula. Thus, the computation of the distribution is done 
in a constant time due to the closed formula while the exact computation
is much longer. Hence, we obtain an approximate solution of the expected
cost very quickly.
In order to find the best threshold $m$ for which it is better to activate
or deactivate, we proceed by comparing approximated costs
of having $k$ VMs activated with the one of having $k+1$ 
(respectively $k-1$)
VMs activated added by the activation (respectively deactivation) cost.

\paragraph{Design of the heuristic}
Let us recall, the formulas of the stationary distributions in $M/M/k/B$ queues. 
Let $\hat{\pi}(k,m)$ be the stationary distribution of state $(m,k)$, we get:  
\begin{equation*}
\hat{\pi}(k,m) = \hat{\pi}(k,0) \cdot \dfrac{\rho^m}{m!} \, \ \mbox{ for } 1 \leq m \leq k
\quad \mbox{ and } \quad 
\hat{\pi}(k,m) = \hat{\pi}(k,0) \cdot \dfrac{a^m \cdot k^k}{k!} \, \ \mbox{ for } k < m \leq B.
\end{equation*}
where $\rho = \dfrac{\lambda}{\mu}$,  $a=\dfrac{\rho}{k} < 1$ and where
$\hat{\pi}(k,0)$ is:
\begin{equation*}
\hat{\pi}(k,0) \! = \!
\Bigg( \sum_{m=0}^{k-1} \dfrac{\rho^m}{m!} + \dfrac{k^k}{k!} (B-k+1) \Bigg)^{\! -1}  \mbox{ if } a = 1 \ \mbox{and} \  
\hat{\pi}(k,0) \! = \!
\Bigg( \sum_{m=0}^{k-1} \dfrac{\rho^m}{m!} + \dfrac{\rho^k}{k!} \dfrac{1 - a^{B-k+1}}{1-a} \Bigg)^{\! -1} \mbox{ if } a \neq 1.
\end{equation*}

Let us define $\hat{C}_k(m)$ as the approximate cost of having $k$ VMs turned on and $m$ requests in the system by:
\[ 
\hat{C}_k(m) =\hat{\pi}(k,m) \bigl( m \cdot C_H + k \cdot C_S \bigr) 
+ \hat{\pi}(k,B) \cdot \lambda \cdot C_R \, ;
\]
define $ \hat{C}^A_{k+1}(m)$ as the approximate cost of having $m$ requests in level $k+1$ knowing that we have just activated a new VM by:
\[
\hat{C}^A_{k+1}(m) = \hat{\pi}(k+1,m) \bigl( m \cdot C_H + (k+1) \cdot C_S \bigr)
+ \hat{\pi}(k+1,B) \cdot \lambda \cdot C_R + \hat{\pi}(k,m-1) \lambda C_A \, .
\]
Let us define $\hat{C}^D_k(m)$ as the approximate cost of having $m$ requests in level $k$ knowing that we have just deactivated a new VM by:
\[
\hat{C}^D_k(m) = \hat{\pi}(k,m) \left( m C_H + k C_S \right) 
+ \hat{\pi}(k,B) \lambda C_R
+ \hat{\pi}(k+1,m+1) (k+1) \mu C_D \, .
\]
We want to compare these costs to know whether we should activate, deactivate or
let the servers unchanged.

For the activation thresholds. We  define 
$\phi^A_k(m) = \hat{C}_k(m) - \hat{C}^A_{k+1}(m)$. 
If $\phi^A_k(m) < 0$ then it is better to not activate a new server 
while if $\phi^A_k(m) \geq 0$ it is better to activate a new one.
There is no evidence of the monotonicity of $\phi^A_k(m)$ in $m$,
nevertheless we choose the threshold by $F_k = \min \left\{ m :  \phi^A_k(m) \geq 0\right\}$.
To compute the whole thresholds, all $k$ are studied in an ascending
order. For a fixed level $k$, we need to compute $\phi^A_k(m)$ for all $m$ and stop when the function $\phi^A_k(m)$ is larger than $0$. 

For deactivation thresholds. We define 
$\phi^D_k(m) = \hat{C}_{k+1}(m) - \hat{C}^D_k(m)$.
If $\phi^D_k(m) < 0$ then it is better to stay with $k+1$ servers while
if $\phi^D_k(m) \geq 0$ it is better to deactivate a virtual machine.
Similarly, there is no evidence about the monotonicity of $\phi^D_k(m)$ but
we define $R_k = \min \left\{ m :  \phi^D_k(m) \geq 0\right\}$.
The computation of the whole deactivation thresholds is similar to the activation ones.


\paragraph{Markov chains approximation methods coupled with local search heuristics}

The initialisation of the local search heuristics presented in \ref{sec:Localsearch} either is totally random or simply chooses
the bounds (lowest values or highest values) of the solution space.
The aim here is to improve the speed and the accuracy of these
heuristics by coupling them with the \textit{MMKB approximation} used as an initialisation step.
This initialisation aims at finding an initial solution that will
be in the basin of attraction of the optimal solution. Starting 
with this solution improves the ability
of local search algorithms to reach the best solution in a faster time.
This new method behaves better than
the initialisation based on \textit{Fluid approximation} adapted from 
\cite{art:mitrani2013managing} in the work \cite{art:Tou19}. 
We define by \textbf{Alg MMK} the coupled heuristic with MMKB approximation.

\section{Computing policies with Markov Decision Process}
\label{sec:MDPmodel}

We consider here the multi-server queue of Section \ref{sec:controlledQueue}  
and the controlled model in which only a single virtual resource can be activated or 
deactivated at a time decision. This optimal control problem is a continuous time process 
and thus should be solved using a Semi Markov Decision Process (SMDP).
In this section, we first describe the uniformised SMDP and its solving.
Then we underline the differences between the optimal control model
and the optimisation problem of Section \ref{sec:heuristicsMCmodel}.

\subsection{The SMDP model}

\subsubsection{Elements of the SMDP}
The state space is the one defined in Section \ref{sec:cloudModel}: $\mathcal{S}=\{0,1,\ldots,B\} \times \{1,\ldots,K\}$.
Hence, a state $x \in \mathcal{S}$ is such that $x=(m,k)$ where  $m$ is 
the number of customers in the system and $k$ the number of active servers.
Similarly, the action space $\A ~ = ~ \{-1,0,1 \}$ represents respectively: deactivate one machine, left unchanged  the active servers, or activate one machine.

The system evolves in continuous time and at some epoch a transition occurs.
When a transition occurs, the controller observes the current 
state and reacts to adapt the resources by activating or deactivating the virtual 
machines. The activation or deactivation are instantaneous and then the system 
evolves until another transition occurs. Two events can occur:  an arrival with 
rate $\lambda$ which increases the number of customers present in the system by 
one or a departure which decreases the number of customers by one.

Let $x=(m,k)$ be the state and $a$ the action, we define $N(k+a)$ 
as the real number of active VM after the action 
$a$ was triggered. We have 
$N(k+a) =  \Min\{\Max\{1,k+a\},K\}$.  It follows that the transition rate is equal to
\begin{equation*}
\begin{cases}
\lambda &  \text{if } y = \! \bigl( \Min\{m  \! + \! 1,B \}, N(k+a)\bigr) \\
\mu \! \cdot \! \min \{m,N(k \! + \!a)\} & \text{if } y = \! \bigl(\Max\{0, m -\! 1\}, N(k+a)\bigr)
\end{cases}.
\end{equation*}

\subsubsection{Uniformised transition probabilities} \label{sec:unif}
We apply here the standard method to deal with continuous time MDP: the 
uniformisation framework. We follow the line of chapter 11 of \cite{liv:put1994} 
to  define the uniformised MDP. 

We define $\Lambda(m,k,a)$ as the rate per state. We have 
$\Lambda(m,k,a) =  \lambda + \mu \! \cdot \! \min \{m, N(k+a) \}$. 
From now on, any component denoted by "$\sim$" refers to the uniformised process.
We thus define the uniformisation rate by $\tilde{\Lambda}$ with 
$\tilde{\Lambda} =  \max_{(m,k,a)}  \Lambda(m,k,a) = \lambda + K \mu$
which is the maximum transition rate.
The transition probability from state $x$ to state $y$ when action $a$ is triggered is denoted by $\tilde{p}(y|x,a)$
in the uniformised model. 

We have:
\begin{equation*}
\tilde{p}(y|x,a)=
\begin{cases}
\dfrac{\lambda}{\tilde{\Lambda}} & \text{if } y=\left(m+1,N(k+a)\right) \\
\dfrac{\mu \min \{m,N(k +a)\}}{\tilde{\Lambda}} & \text{if } y=\left(m-1,N(k+a)\right) \\
\dfrac{\tilde{\Lambda}-\Lambda(m,k,a)}{\tilde{\Lambda}} & \text{when } y=(m,k) \\
0 & \text{otherwise} \, ,
\end{cases} 
\end{equation*}
when $x=(m,k)$ such that $0 < m < B$ and $a$ is arbitrary ; 
when $x=(B,k)$ with $a \not = 0$, or 
when $x=(B,k)$ with $k \not =K$ and $a \not = +1$ 
or also when $x=(B,k)$ with $k \not =1$ and $a \not = -1$; 
when $x=(1,k)$ with $a \not = 0$, or  
when $x=(1,k)$ with $k \not =K$ and $a \not = +1$, or also  
when $x=(1,k)$ with $k \not =1$ and $a \not = -1$. 
We have: 
\begin{equation*}
\tilde{\Lambda} \times \tilde{p}(y|x,a)=
\begin{cases}
\mu \min \{m,N(k +a)\} & \text{if } y=\left(B-1,N(k+a)\right) \\
\tilde{\Lambda}-\mu \min \{m,N(k +a)\} & \text{when } y=(B,k) \\
0 & \text{otherwise} \, ,
\end{cases} 
\end{equation*}
when $x=(B,k)$ with $a=0$, or 
when $x=(B,k)$ with $k \not =K$ and $a \not = +1$, 
or also when $x=(B,k)$  with $k \not =1$ and $a \not = -1$.
We have
\begin{equation*}
\tilde{\Lambda} \times \tilde{p}(y|x,a)=
\begin{cases}
\lambda & \text{if } y=\left(1,N(k+a)\right) \\
\tilde{\Lambda}-\lambda & \text{when } y=(0,k) \\
0 & \text{otherwise} \, ,
\end{cases} 
\end{equation*}
when $x=(0,k)$ with $a=0$, 
or when $x=(0,k)$ with $k \not =K$ and $a \not = +1$, or  
when $x=(0,k)$  with $k \not =1$ and $a \not = -1$.

\subsubsection{The uniformised stage costs}
We take the definition of costs of Section \ref{sec:controlledQueue}.
Instantaneous costs are charged only once and
are related to activation, deactivation and losses.
Accumulated costs are accumulated over time 
and are related to consumption and holding cost.
After uniformisation of instantaneous and accumulated costs (Chapter 11.5 of \cite{liv:put1994}), we obtain the following equation for $x=(m,k)$:
\[
\tilde{c}(x,a)= 
\bigl( C_A  \mathds{1}_{\{a=1\}} + C_D  \mathds{1}_{\{a=-1\}}\bigr)   \dfrac{\Lambda(x,a)}{\tilde{\Lambda}} 
+ \dfrac{\lambda}{\tilde{\Lambda}} C_R \mathds{1}_{\{m = B \}}
+ \dfrac{1}{\tilde{\Lambda}} 
\bigl( N(k+a) \cdot C_S  + m \cdot C_H \bigr) \, .
\]

\paragraph{Objective function}
The objective function defined in Equation \eqref{eq:globalobjective}
translates for all $ x \in \mathcal{S}$  into 
\begin{equation*}
\rho^{\pi}(x)  =  \lim_{N \rightarrow \infty} E^{\pi} 
\Bigl[ \dfrac{1}{N} \sum \limits_{t=0}^{N-1} \tilde{c}(x_t,\pi(x_t)) ~ | ~ x_0  =  x \Bigr] \, ,
\end{equation*} 
for a  given policy $\pi$. The value $\rho^{\pi}$ is the expected stage
cost and equals $\bar{C}$ when actions follow policy $\pi$.

\subsection{Solving the MDP}
\subsubsection{Classification of the SMDP}
Our SMDP is an average cost model. This is why, the expected stage costs depend on the recurrent properties of the underlying Markov chain that is generated by a deterministic policy. A classification is then necessary to study them. We use the classification scheme of  \cite{liv:put1994}.

\begin{definition}[Chapter 8.3.1 in \cite{liv:put1994}]
A MDP is:
\begin{enumerate}
\item[$i$)] \emph{Unichain} if, the transition matrix corresponding to every deterministic stationary policy is unichain, that is, it consists of a single recurrent class, plus a possibly empty set of transient states;
\item[$ii$)] \emph{Multichain} if, the transition matrix corresponding to at least one deterministic stationary policy contains two or more recurrent classes;
\item[$iii$)] \emph{Communicating} if, for every pair of state $x$ and $y$ in ${\cal S}$, there exists a deterministic stationary policy $\pi$ under which $y$ is accessible from $x$, that is, $p_{\pi}^n(y,x) >0$ for some $n \geq 1$;
\end{enumerate}
\end{definition}

\begin{proposition}
\label{prop:multichain}
The MDP is multichain. There is a stationary deterministic 
policy with monotone hysteresis properties that induces a corresponding  Markov chain with more than two different recurrent classes.
\end{proposition}
\begin{proof}
We assume that $K \geq 2$ and let $k$ be such that $k \in [1, \ldots, K] $. We define
the policy $q$ as follows.
For any level $l$ with $l < k-1$, we have only activation.
There exists $m \in [0, \cdots, B]$ such that $q(m',l) = 1$ for any 
$m \leq m'$ and $q(m,l)\leq q(m,l-1)$.
For level $k$ we have neither activation nor deactivation and  
$q(m,k) = 0$ for all $m \in [0, \cdots, B]$.
For any level $l$ with $l>k+1$, we have only deactivation.
There exists $m \in [0, \cdots, B]$ such that $q(m',l) = -1$  for any $m' \leq m$
and  $q(m,l)\geq q(m,l+1)$. Therefore, we have three recurrent classes: the level $k$, the level $k-1$ and the level $k+1$ (and two recurrent classes for $K=2$). 
See supplementary materials in appendix \ref{app:analysis}.
\end{proof}

\begin{lem}
\label{prop:communicating}
The MDP is communicating. There exists a stationary isotone hysteresis 
policy such that the corresponding Markov chain is irreducible.
\end{lem}
\begin{proof}
We exhibit such a policy. Let $\pi$ be the deterministic stationary policy such that the thresholds $l_k$ are defined by $l_k=0$ and the thresholds $L_k$ by $L_k=B$ for all $k$. The induced Markov chain is irreducible since any level can be reached from another one (when $m=0$ or $m=B$) and since in a given level all the states are reachable from any state. Thus, we have that $p_{\pi}^n(y,x) >0$ for some $n \geq 1$ for all couples $(x,y)$. Therefore the MDP is communicating.
\end{proof}

\paragraph{Bellman Equations} 
In the \emph{multichain} case, the Bellman Equations are composed by two equations.
In the uniformised model, we then have (\cite{liv:put1994}) the two following optimality equations:
\begin{equation*} 
\min_{a \in A}  \Bigg\{
 \sum \limits_{y \in \X} \tilde{p}(y \, | \, x,a) \rho(y) - \rho(x)
\Bigg\} =0 
\ \text{ and } \ 
U(x) =  \min_{a \in B_x} \Bigg\{ \tilde{c}(x,a) - \dfrac{\rho(x)}{\tilde{\Lambda}}  +  \sum \limits_{y \in \X} 
\tilde{p}(y \, | \, x,a) \cdot U(y) \Bigg\} 
\end{equation*}
for all $x \in \X $, where 
\begin{equation*}
B_x  =  \Big\{ a \in \A ~ \big| ~  \sum \limits_{y \in \X} 
\tilde{p}(y \, | \, x,a) \rho(y)  ~ = ~ \rho(x) \Big\} \, .
\end{equation*}
It could be noticed that in the unichain case, $B_x ~ = ~ \A$ and that the two equations reduce to only the second one. 
These two non linear equation systems should be numerically solved 
to find $U(x)$ and $\rho(x)$. Once these terms are approximated we deduce the optimal policy with: 
\begin{equation*}
q(x)= \arg \min_{a \in B_x}  \bigg\{ \tilde{c}(x,a)
+\sum \limits_{y \in \X} \tilde{p}(y \, | \, x,a)\cdot U(y)  \bigg\}  \,.
\end{equation*}

\subsubsection{Algorithms}
\label{sec:MDPalgos}

We describe here the choice of the algorithms used to solve this multichain and communicating model. Computing multichain model is much complicated namely since
testing if an induced Markov chain is unichain is a NP complete problem.
We first show that we actually can use some unichain algorithms due to the communicating properties of our SMDP, then we present two structured algorithms based on hysteresis properties.

\paragraph{Unichain algorithms}
From \cite{liv:put1994} it exists a \emph{multichain policy iteration algorithm} that solves multichain models. It requires to solve two equation systems and thus is time consuming. However, by Proposition \ref{prop:communicating}, our MDP is communicating
and in Theorem 8.3.2 of \cite{liv:put1994} it is proved that  unichain
value iteration algorithm converges to the optimal value in communicating models.
This property allows us to use the unichain value iteration algorithm. 

There is also, in \cite{liv:put1994}, a policy iteration algorithm  for communicating models. It requires to start with an initial unichain policy and its inner loop roughly differs from the unichain case in order to keep some properties.  
We decide to use here algorithms based on unichain policy iteration. There does not exist theoretical guarantee of their convergence, but we showed in numerical experiments that they always converge to the same policy than value iteration.

Four different usual unichain algorithms will be considered: \textit{Value Iteration}, \textit{Relative value iteration}, \textit{Policy Iteration modified} and \textit{Policy Iteration modified adapted}, which adapts its precision in the policy evaluation step.  They are respectively referred by \textit{VI}, \textit{RVI}, \textit{PI} and \textit{PI Adapt}. 
The algorithms are all described in \cite{liv:put1994} and  are already implemented in the software \cite{conf:ajmMarmote}.

\paragraph{Structured Policies algorithm}
We now integrate hysteresis properties in the algorithms.
Two classes of policies have been investigated: Double Level class (Definition \ref{def:doubleLevel}) and Monotone Hysteresis class (Definition \ref{def:hysteresisSVM}). 
The goal is to plug hysteresis assumptions during the policy improvement
step of policy iteration (Policy Iteration has two major steps). This allows to test 
less actions at each iteration and to speed up the algorithm.
Two algorithms are implemented: one for the double level properties 
(referred as \textit{DL-PI}) and one for the hysteresis properties
(referred as \textit{Hy-PI}). 

On the other hand, we do not have theoretical guarantee that these methods converge, first because we do not theoretically know if 
hysteresis policies are optimal for our model, and second because the underlying
PI also has no convergence guarantees in multichain.
Nevertheless, in the numerical experiments we made (see later) 
all the optimal policies returned by classical algorithms 
have hysteresis properties. Furthermore, all MDP algorithms considered here (structured as well as classical) returned the same solution.
This therefore underlines the interest of considering such hysteresis policies 
especially since the gain in running time is obvious as observed in the experiments.

\paragraph{Computation of the hysteresis thresholds}
The non-structured MDPs (also called simple MDPs) do not assume any restrictions for their policy research. Thus, they return the optimal policy and so return a decision rule  $q^*$ which gives the optimal action to take but not the
thresholds. Therefore, we need to test the hysteresis property and to
compute the $l$ and $L$ hysteresis thresholds consistently with Definition \ref{def:hysteresisSVM}.
Let $q^*$ be the optimal policy returned by the PI algorithm.
We check if  $q^*$ is monotone, if not the optimal policy is not hysteresis.
If so, we proceed as follows. We consider, for all $k$, the set 
$\{ m \mid q^*(m,k) = 1 \text{ and } q^*(m-1,k) = 0 \}$. If the set is of size $2$ then
the policy is not hysteresis else it is hysteresis. 
In hysteresis case, if the set is of size $1$ then 
$L_{k+1} =m$ and if the set is empty then $L_{k+1} =  \infty$.
Also, we consider, for all $k$, the set 
$\{ m \mid q^*(m+1,k) = 0 \text{ and } q^*(m,k) = -1 \}$.
If the set is of size $2$ then
the policy is not hysteresis, else it is. If the set is of size $1$ then 
$l_{k+1} =m$ and if the set is empty then $l_{k+1} = 0$.

\subsection{Theoretical comparison between the two approaches MC and MDP}
\label{sec:modelDifferences}
Although they seem very similar these two models present some rather subtle theoretical differences with notable consequences. They are studied here while the numerical comparisons will be carried out in Section \ref{sec:numericalExp}.

\subsubsection{Number of activated resources}
\label{sec:MDPresourcemanagement}
A major difference between the two approaches is related to the 
class of policy they consider and the consequences on the management of resources. 
The (MC) approach (Section \ref{sec:MCmodel}), deals with 
\emph{strictly isotone} policies 
(\textit{i.e.}
$0= l_1 < l_2 < \ldots < l_{K} \leq B$   and  
$0 = L_1 < L_2 < \ldots < L_{K} \leq B$).
Therefore, $L_K$ can not be infinite, and in this model all the  $K$
servers should be activated. Furthermore, since inequalities are strict 
$l_K \geq K$, the only state in which there is only one active server 
is restricted to the level $k=1$.
Strictly isotone assumption is required to keep the size of the chain constant
(see \cite{art:ibekeilson95}). 
On the other hand, the (MDP) approach either does not assume any 
structural property on the policy or considers isotone policy (see Definition \ref{def:isotoneHysteresis}). In this approach constraints are looser and 
it is possible not to have all the servers activated and conversely to have many resources activated even if the queue is empty. 

The resource management follows the optimal policies returned by each of the approaches. Three categories of resource management have been identified.  
We noticed that these categories depend on the relative values of the parameters ($\lambda, \mu$) between them. However, we do not know how to precisely quantify their borders. 
\begin{definition}
\label{def:MDPscenarios} These categories are:
\begin{enumerate}
    \item \emph{Medium arrival case}: All {\it VMs} are turned On and turned Off 
    in both approaches.
    This occurs when the system load is \textit{medium}, i.e. the arrival rate $\lambda$ is close to the service rate $k \times \mu$.  
    \item \emph{Low arrival case}:  All {\it VMs} are turned On and turned Off in (MC) while in (MDP) some {\it VMs} will never be activated.
     This occurs when the system load is \textit{low}, i.e. the arrival rate $\lambda$ is very small compared to the service rate $k \times \mu$.
    \item \emph{High arrival case}:  All {\it VMs} 
    are turned On and turned Off
    in (MC) while in  (MDP) some {\it VMs} will never be deactivated.
    This occurs when the system load is \textit{high}, i.e. the arrival rate $\lambda$ is very large compared to the service rate $k \times \mu$.
\end{enumerate}
\end{definition}

\begin{proposition}[Non optimality of strictly isotones policies]
The above classification helps us to claim that strictly isotone policies are not optimal and that MDP approach performs better. 
\end{proposition}
\begin{proof}
During the numerical experiments in Section \ref{sec:numericalExp},
we identify numerous cases in which strictly isotone policies are not optimal. This is mainly due to the phenomenon of non activation or non deactivation exposed above: examples of non-optimality are found in both  low arrival and high arrival categories.
The structured (MDP)s have less constraints and return the optimal solution which is the one computed by simple MDPs.
See appendix \ref{app:example} for details.
\end{proof}

The proposition above highlights another benefit of MDP approaches since
they allow to size the exact number of machines to be activated (this is
particularly true in the \emph{low arrival} case). Hence, if in a policy
resulting from a MDP there exits
a $k$ such that for all $l >k $ and for all $m \in \{1,\ldots,B\}$ 
we have $q(m,l)<1$, then $K-k$ machines are not necessary.
This is not true in MC heuristics,
nevertheless, it is possible to adapt the previous heuristics to find the optimal number of VM to be activated in \emph{low arrival} case.  This requires running the algorithms for each level $k$ with $k \in \{1,\ldots,K\}$. Then take
the level $k$ which has the smallest cost says $k^*$. If $k^* =K$ all machines 
must be turned on in $K$, otherwise only  $k^*$ machines are important and the other ones useless. But this method is time consuming. 
Also, such an adaptation for \emph{high arrival} cases is not obvious.

\subsubsection{Different temporal sequence of transitions}
\label{temp_diff}
The last point to investigate is the difference in dynamical behavior between transitions. This difference is slight and has no effect on average costs, nevertheless it induces a difference on the values of the thresholds.
In this system, there are two kinds of transitions: 
\emph{natural transitions} which are due to events (departures or arrivals) and
\emph{triggered transitions} which are caused by the operator (activation or
deactivation). These two transitions are instantaneous.
Due to their intrinsic definitions, the models considered here do not
observe the system at the same epochs.  Hence, in the Markov chain the system is observed just before a natural transition while in the MDP the system is observed just after a natural transition. More formally, let us assume that $x$ is the state before a natural transition: it is seen by the Markov chain. Then the transition occurs and the state changes instantaneously in $x'$ that is the state seen by the MDP. The controller reacts and the triggered transition occurs, thus the system moves
in $x''$ in which the system remains until the next event which will cause the next
natural transition. Since state changes are instantaneous, this has no impact on
costs. 
See appendix \ref{app:analysis} for details.

In the \textit{medium case arrival} defined by Definition \ref{def:MDPscenarios}, thresholds are fully comparable and we have the following lemma
\begin{lem}
In the \textit{medium case arrival}, hysteresis thresholds of Markov chain and SMDP hysteresis threshold can be inferred from each other.
Let  $(F_1,\ldots,F_{K-1})$ and $(R_1,\ldots,R_{K-1})$ be the thresholds of the MC model and let $(L_2,\ldots,L_{K})$ and $(l_2,\ldots,l_{K})$ be the SMDP thresholds.
Then, for all $k \in \{1,\ldots,K-1\}$, we have: $L_{k+1}-1=F_k \quad \text{ and } \quad l_{k+1}+1=R_k.$
\end{lem}
\begin{proof}
We make the proof for an activation threshold. 
Let assume that the state just before a transition is $x=(F_k,k)$, if an arrival event occurs then the system moves instantaneously in a state $x'=(F_k+1,k)$.
According Definition \ref{def:hystGolub}, a virtual machine is activated and then 
the system instantaneously moves again in $x''=(F_k+1,k+1)$. Observe now, that the
state $x'$ (according Definition \ref{def:monotoneHysteresis}) 
is the state in which  the SMDP
observes the system an takes its activation decision. Thus, we have $L_{k+1}=F_k+1$.
The proof works similarly for deactivation. 
\end{proof}

\section{Real model for a Cloud provider}
\label{sec:concreteModel}
Now we are going to give real environment values to our models.
The approach is difficult since, up to our knowledge, 
it does not exist,
in the literature, an unified model including both energy, quality
of service, and real traffic parameters. 
We propose to build a global cost taking into account both
real energy consumption, real financial costs of virtual machines
and Service Level Agreement (SLA) simultaneously 
from separate works.
We think that the model presented is sufficiently generic
to represent a large class of problems and is enough to give
a relevant meaning to all parameters.
With our algorithms and this model 
the cloud owner can generate meaningful optimised costs in real cases. Parameter values used in experiments can be found in Section \ref{sec:ResNumReel}.

\subsection{Cost-Aware Model with energy consumption traces and real prices}

\paragraph{Measurements of energy consumption}
Many measurements are performed  on the real datacenter 
\textit{grid5000} in \cite{art:benoit}, and 
energy consumption data are obtained for VMs hosted 
on physical servers in \cite{art:orgerie}.
We keep these values in Watt for energy
consumption of virtual machines. 

\paragraph{Financial Cost}
In order to keep the relevance of our \emph{cost-aware} model,  we must represent a financial cost. Henceforth, we transform energy consumption given in  Watts  into a financial cost based on the price in euros of the
KWh in France fixed by national company. 
To obtain the financial values of the other prices, we observe the commercial offers of providers \cite{rap:aws2}. This gives us 
the operational costs. 

\paragraph{Service Level Agreement}
We propose to modify the performance part of our cost function 
for capturing the realistic scheme of the Service Level Agreement. 
Actually, we have to give a concrete meaning to our
holding cost since it does not translate directly from usual SLA models. 
Indeed, in SLA contracts, it could be stated that
when the response time exceeds a pre-defined
threshold $T_{SLA}$ then penalty costs must be paid. 
Here, the penalty $C_p$ is the price that a customer pays for one hour of
service of a virtual machine (full reimbursement) \cite{rap:aws2} and
the threshold corresponds to a maximal time to process a request.  

With our model, it is possible to define an holding cost that
models the penalty to pay when the QoS is not satisfied.
We focus on the mean response time in the system  
where the response time is denoted by $R_T$.
The SLA condition translates in:  $E[R_T] \leq T_{SLA}$. We introduce  the mean number of customers $E(N)$ and with Little's law we get:
$E[N] \leq T_{SLA} \cdot \lambda$.
We thus may define a customer threshold $N_{SLA}$ such that
$N_{SLA} ~ = ~ T_{SLA} \cdot \lambda$.
This customer threshold is the maximum number of requests, 
that the system can accept to keep 
the required \textit{QoS} satisfied.
Therefore, each time $m > N_{SLA}$, the cloud provider will pay a holding cost penalty of $C_p$ per customer.

\paragraph{Global cost formula}
We obtain the expected global cost charged to the cloud provider:
\begin{align*}
\overline{C} (m , k)&  = C_p \cdot (m - N_{SLA},0)^+  
+ C_S \cdot k 
+ C_p \cdot \lambda \cdot \mathds{1}_{\{m=B,k=K\}} 
+ C_A \cdot \lambda \cdot \mathds{1}_{\{m = F_k,\, k < K\}} \\
& \qquad + \: C_D \cdot \mu \cdot 
\min\{m,k\} \cdot \mathds{1}_{\{m = R_{k-1}+ 1, \; 2 \leq k \leq K \}} + C_{static} \, ,
\end{align*}
In the energy consumption part, $C_{static}$ represents the static financial cost derived from the idle energy consumption of the physical server which hosts the virtual machines.  
The remaining energy part of the cost function is left unchanged with the activation, deactivation, and energy costs of a VM defined above.

\subsection{Real packet traffic  and CPU utilisation}
For obtaining concrete values of parameters $\lambda$ and $\mu$ we must search workload traces with real scenarios in public cloud.  
Here the traces come from MetaCentrum Czech National Grid data \cite{art:orgerieWL}.

\section{Numerical experiments}
\label{sec:numericalExp}
This part is devoted to numerical experiments and the comparison between all previous algorithms. We have implemented all methods in C++, with the stand-alone library \cite{conf:ajmMarmote}.

\subsection{Experiments design and results}

All the experiments were carried out on a platform 
built on an Intel Xeon
X7460 processor, with $6 \times 4$ cores at $2,66$GHz and $16$
GB of RAM. 
Numerical experiments have been done for different Cloud model parameters. System parameters were taken arbitrarily. We defined four scenarios with different scales: \textbf{Scenario A}: K=3, B=20;
\textbf{Scenario B}: K=5, B=40; \textbf{Scenario C}: K=8, B=60;
\textbf{Scenario D}: K=16, B=100.
We call \textit{Instance} a set of parameters $\lambda, \mu, C_a, C_d, C_h, C_s, C_r$. The algorithms have been launched
on $46656$ instances.
Costs ($C_a, C_d, C_h, C_s$) were taking values in $[0.5, 1, 2, 5, 10, 20 ]$, $C_r$ in $[1, 10, 100, 1000, 5000, 10000]$, queueing parameters $\lambda, \mu$ in $[0.5, 1, 2, 5, 10, 20]$. 
We rank the different algorithms by assessing 
their accuracy and their running time.
We first studied numerical experiments in the (MC) approach by comparing  the different heuristics by these two criteria. 
For Scenario A, accuracy is computed based on optimal solution generated by exhaustive search. For higher scale scenarios (B,C,D), we compute accuracy regarding the best solution among the heuristics.
Then, we compared all the (MDP) algorithms with the value iteration solution (convergence is known from Section \ref{sec:MDPalgos}) 
before comparison with the heuristics to compare both approaches. 
We finally ran the experiments for a concrete cloud scenario with
real data, showing how we can effectively benefit from the best
methods to calculate optimal thresholds. 

\begin{table}[hbtp]
\resizebox{0.75\textwidth}{!}{ %
\begin{minipage}{\textwidth}
\hspace{0.5cm}
  \begin{tabular}{c|c|cc||cc||cc||cc}
   \textbf{Models} &
   \textbf{Algorithms} &
      \multicolumn{2}{c||}{\textbf{Scenario A}} &
      \multicolumn{2}{c||}{\textbf{Scenario B}}  &
      \multicolumn{2}{c||}{\textbf{Scenario C}} & 
      \multicolumn{2}{c}{\textbf{Scenario D}} \\
      \hline
      & & Time (sec) & \% Opt  & Time (sec) & \% Min & Time (sec) & \% Min  & Time (sec) & \% Min \\
      \multirow{5}{*}{MC Heuristics} 
    & \textbf{BPL} & 0,063 sec & 96,9 \%  & 2,684 sec & 86,81\% & 37,23 sec & 79,57 \% & 757 sec & 55,39\% \\ 
    & \textbf{BPL Agg} & 0,047 sec & 96,9 \% & 1,163 sec & 86,81\% & 10,27 sec & 79,57 \% & 186 sec & 55,39\% \\ 
    & \textbf{{BPL-MMK Agg}} & 0,036 sec & 97,78 \% & 0,53 sec &  96,59 \%  & 5,16 sec & 95,39 \% & 31,7 sec & 87,24 \% \\
    & \textbf{NLS} & 0.043 sec & 92,70 \% & 1,007 sec & 62,24 \% & 12,38 sec & 43,51 \% & 257 sec & 20,44 \%  \\
    & \textbf{NLS Agg} & 0,034 sec & 92,70 \% & 0,458 sec & 62,24 \% & 4,05 sec & 43,51 \% & 61 sec & 20,44 \%  \\
    & \textbf{{NLS-MMK Agg}} & 0,014 sec &  96,32 \% & 0,06 sec & 93,58 \%  & 0,9 sec & 89,54 \%  & 8,39 sec & 76,17 \% \\
    \hline
    \hline
    
    \multirow{7}{*}{MDP} 
    
    & & Time (sec) & \% Opt  & Time (sec) & \% Opt & Time (sec) & \% Opt  & Time (sec) & \% Opt \\
    
    & \textbf{VI} & 0.0057 sec & 100 \%  & 0.021 sec & 100 \% & 0.0406 sec & 100 \% & 0,0944 sec & 100 \% \\ 
    & \textbf{RVI} & 0.0057 sec & 100 \% & 0.021 sec & 100 \% & 0.0406 sec & 100 \%  & 0,095 sec & 100 \% \\
    & \textbf{PI} & 0.0028 sec & 100 \% & 0.0124 sec & 100 \% & 0.0214 sec & 100 \%  & 0.0606 sec & 100 \%\\
    & \textbf{PI Adapted} & 0,0023 sec & 100 \% &  0,0117 sec & 100 \% & 0,0201 sec & 100 \% & 0,0583 sec & 100 \% \\
    & \textbf{DL-PI} & 0,00115 sec & 100 \%  & 0,0072 sec & 100 \% & 0,0105 sec  & 100 \% & 0,0452 sec & 100 \% \\
    & \textbf{Hy-PI} & 0,00113 sec & 100 \%  & 0,0069 sec & 100 \% & 0,0100 sec & 100 \% & 0,0442 sec & 100 \% \\

    \hline
    \hline

    \multirow{5}{*}{Comparison}
    & & Time (sec) & \% Opt  & Time (sec) & \% Opt & Time (sec) & \% Opt  & Time (sec) & \% Opt \\
    & \textbf{{NLS-MMK Agg}} & 0,014 sec & 64,82 \%  & 0,06 sec & 61,52 \% & 0,9 sec & 59,82 \% & 8,39 sec & 52,61 \% \\
    & \textbf{{BPL-MMK Agg}} & 0,036 sec & 65,15 \% & 0,53 sec & 62,65 \% & 5,16 sec & 61,69 \% & 31,7 sec & 57,1 \%\\
    & \textbf{DL-PI} & 0,00115 sec & 100 \%  & 0,0072 sec & 100 \% & 0,0105 sec  & 100 \% & 0,0452 sec & 100 \% \\
    & \textbf{Hy-PI} & 0,00113 sec & 100 \%  & 0,0069 sec & 100 \% & 0,0100 sec & 100 \% & 0,0442 sec & 100 \% \\
  \end{tabular}
  \end{minipage}}
  \caption{Numerical experiments for comparison of heuristics and MDP algorithm}
  \label{tab:Results}
\end{table}

\subsection{Experiments analysis}

We display in Table \ref{tab:Results} the comparisons for (MC) heuristics and (MDP) algorithms. For each algorithm, we display for each scenario the running time in seconds in the first column and the accuracy in the second column.

\subsubsection{Assessing heuristics approaches}
We want to assess the selected heuristics of \cite{art:Tou19} and
the gain of the improvements proposed.

\paragraph{Aggregation strongly improves algorithms time speed}
For the same efficiency, we first observe a significant time saving provided by the aggregation technique. As expected, the improvement in Scenario A is small, but the gain increases as the size of the system rises. Henceforth the execution time is divided by about $1.3$ when $K=3$ and $B=20$ and by about $3.8$ when $K=16$ and $B=100$ (e.g. BPL algorithm passes from 757 seconds to 186 seconds).


\paragraph{Impact of the M/M/k/B approximation}
Coupled with heuristics \textbf{BPL} and \textbf{NLS}, M/M/k/B approximation brings significant improvement considering the efficiency and the time execution. For a large scale case (Scenario D in Table \ref{tab:Results}), the \textit{BPL MMK Agg} provides the best solution for 87.24~\% of the instances in 31.7 seconds, which is the best heuristic in terms of time-accuracy ratio. We can notice, in all cases, that coupled heuristic always have better accuracy than the heuristic alone. Moreover, the initialisation technique does not only improve efficiency but also time execution. For example in Scenario C, \textbf{BPL Agg} has a mean time of 10,27 seconds for 79,57 \% accuracy while \textbf{BPL-MMk Agg}
obtains 95,39 \% accuracy in 5,16 seconds, 
providing a double major gain.

\subsubsection{Comparison between MDP algorithms}

By Section \ref{sec:MDPalgos}, we theoretically know that value iteration algorithm converges to the exact solution. 
We observe in these numerical experiments, that all MDP algorithms obtain $100$\% of accuracy (optimal solution being given by value iteration) in all scenarios. Thus, all MDP algorithms converge to the optimal solution. 
Concerning the running time, we observe that  policy iteration algorithms \textbf{PI} and \textbf{PI Adapted} are twice as good than value iteration on all studied scenarios.
Furthermore, we see that the execution time of the structured MDP algorithms \textbf{DL-PI} and \textbf{Hy-PI} is 
divided by $2$ compared to \textbf{PI Adapted}.
This leaves us to conclude that integration of policy with structural properties strongly accelerates the convergence. 


\subsubsection{Comparison between (MDP) and (MC) approaches}
We now compare which approach works best to obtain optimal
thresholds.
We first observe that the accuracy difference is
striking, indeed, the MC approach is optimal only about
half the time. Indeed, from Section
\ref{sec:MDPresourcemanagement}, we know that
depending on the queueing parameters, (MDP) approach
could not activate or deactivate certain virtual
resources whereas the (MC) approach always find
thresholds for all levels due to constraints on the
searched policy. This is why (MC) heuristics are
outperformed by (MDP) algorithms.
Hence, for many instances their optimal cost is higher
due to additional activation or deactivation thresholds,
generating extra costs. 
We see that MDP algorithms are faster than heuristics in
the four scenarios. Hence for small scenarios
\textbf{Hy-PI} is about $100$ times faster than the best
heuristic and about $200$ times faster for large
scenarios. Therefore, the MDP approach is
significantly better than the (MC) approach.

\subsubsection{High-scale simulations MDP}

We ran numerical simulations for high-scale scenarios (real size datacenter) with best MDP algorithm \textbf{Hy-PI} to assess its
performance. Simulations were performed on pre-selected costs and
queuing parameters. 
Results are displayed in Table \ref{tab:high-scale}. They show that for small data centers ($64$ VM and $400$ customers)
the optimal policy can be computed in a reasonable time in practice: $2$ seconds. However, for large data centers some further researches should be done since the time represents around $6$ hours. Note that our sizes $K$ and $B$ are larger than in previous works \cite{art:mitrani2013managing}.

\begin{table}[H]
\centering
\resizebox{1\textwidth}{!}{
    \begin{tabular}{p{0.1\textwidth}|p{0.1\textwidth}|p{0.1\textwidth}|p{0.1\textwidth}|p{0.1\textwidth}|p{0.1\textwidth}|p{0.1\textwidth}|p{0.1\textwidth}|p{0.1\textwidth}|p{0.1\textwidth}|p{0.1\textwidth}|}
        & K=3, B=20 & K=5, B=40 & K=8, B=60 & K=16, B=100 & K=32, B=200 & K=64, B=400 & K=128, B=800 & K=256, B=1600 & K=512, B=3200 & K=1024, B=6400 \\ \hline
        Execution time (sec) &  0,00021 & 0.00188 & 0.02739 & 0.05167 & 0.30225 & 2,42377 & 37,97133 & 329.54122 & 2327 & 22850 \\
    \end{tabular}}
    \caption{High-scale datacenter resolution with MDP algorithm \textbf{Hy-PI}}
    \label{tab:high-scale}
\end{table}

\subsection{Numerical experiments for concrete scenarios}
\label{sec:ResNumReel}

\paragraph{Data and scenarios}
We define now numerical values to the real model of Section \ref{sec:concreteModel}. Several case-studies as well as several performance and metrics curves are used to select the values. 
Throughout this section the time unit is the hour. 
The platform comes from \cite{art:orgerie}: the physical hardware
is a Taurus Dell PowerEdge R720 server having a processor Intel Xeon e5-2630 with $6\times 2$ cores up to $2.3$GHz. 
The processor hosts virtual machines by an hypervisor. One core is considered to be one vCPU and cores are distributed equally among the virtual machines:
$12$ VMs with $1$ core per VM, 
$6$ VMs with $2$ cores per VM, and
$3$ VMs with $4$ cores per VM. These three different cases are selected such that we can use the prices of Amazon EC2 instances \cite{rap:aws2}. 

Energy consumption values come from \cite{art:orgerie}.
Physical hardware consumption is around $100$W.
The dynamic consumption of the virtual machines is stated 
in accordance with the number of cores per VM.  Namely, the Watt consumption for activation, deactivation and use is different as it depends on the characteristics of the virtual machines hosted on the physical server. Hence, $1$ VM with $4$ vCPU which executes a request requires $40$W to work while its activation and deactivation take a power of $10$W. At last, all the energy consumption values are transformed in financial costs.

Some workload samples which come from real traces of the MetaCentrum Czech national grid are detailed in \cite{art:orgerieWL}. We select some samples to build concrete daily scenarios for huge size requests.
Hence, we pick a number of arrivals per hour of $\lambda = 50$ and number of served requests per hour as $\mu = 5,10$ or $20$.
Table \ref{tab:scenario} details all the parameter values for each model.

\begin{table}[hbtp]
\small
\caption{Parameter values for three scenarios}
\centering
\begingroup
\setlength{\tabcolsep}{10pt} 
\renewcommand{\arraystretch}{1} 
\resizebox{8cm}{!}{ %
\begin{tabular}{cccc}
\toprule
  Parameter & Model A & Model B & Model C  \\ \midrule
  $\# VMs$ & 3 &  6 & 12  \\ 
  Instances & a1.xlarge & t3.small & t2.micro \\
  $\# vCPUs $ & 4 & 2 & 1 \\
  RAM(Go) & 8 & 2 & 1 \\ \midrule
  $C_p$ & 0.0914\euro /h & 0.0211\euro /h & 0.0118\euro /h \\
  $C_S$ & 0.00632\euro /h & 0.00316\euro /h & 0.00158\euro /h \\
  $C_A = C_D$ & 0.00158\euro  & 0.00079\euro\ & 0.00032\euro \\ 
  $C_{Static}$ & 0.0158\euro /h & 0.0158\euro /h & 0.0158\euro /h \\ \midrule
  $B$ &  $100$ & $100$ & $100$  \\
  $\lambda$ & $50$ req/h & $50$ req/h & $50$ req/h\\
  $\mu$ & $20$ req/h & $10$ req/h & $5$ req/h \\
  \bottomrule
\end{tabular}
}
\endgroup
\label{tab:scenario}
\end{table}

We provide the mean run time of the best algorithm (\textbf{Hy-PI}) among concrete experiments: 0,00127 s for Model A ; 0,0073 sec for Model B and 0,0372 sec for Model C. The optimal solution is computed in a short-time and therefore can allow to recompute new policy when the demand is evolving, in order to dynamically adapt to the need.

\paragraph{Mean global costs given SLA parameter or arrival rate}
From numerical simulations it could be noticed first that, for a fixed arrival rate $\lambda$, when the SLA is not so strict,
then the cost decreases. Second that, when the response time set by the SLA rises, then the operating costs for the cloud providers are reduced. Indeed, providers will pay fewer penalties to clients and fewer power bills for running VMs. On the other hand, we observe that Model C is always better than Model B and Model A for any possible values of $N_{SLA}$ and for any values of $\lambda$ given a fixed SLA parameter. We can conclude that it is always better to decompose physical machine in several virtual resources, since it allows a more efficient resource allocation in dynamic scenario. Indeed, the more virtual resources you have, the less you will pay per utilisation since they require less CPU consumption. Moreover, it allows a more flexible system where you can easily adapt virtual resources to the demand.
See Appendix \ref{app:experiments} for details and more experiments.

\paragraph{Thresholds given SLA parameter or arrival rate}
We also observe that Model C always activates a second virtual resource before the two other models B and A, for any fixed values of $\lambda$ and SLA parameter $N_{SLA}$.
Furthermore, when the arrival rate $\lambda$ is fixed, then all models will activate later when the SLA is less restrictive.
Lastly, when the SLA parameter is fixed, we notice that the first activation threshold decreases when $\lambda$ increases. Indeed, if the demand is growing, it requires more resources and requires faster activation of virtual resources.

\begin{figure}[hbtp]
\centering
\begin{multicols}{2}
    \includegraphics[width=0.33\textwidth]{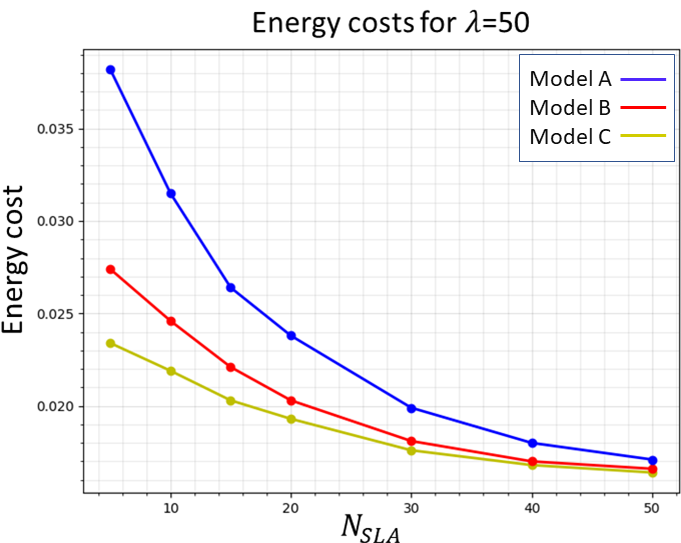}\par
    \includegraphics[width=0.33\textwidth]{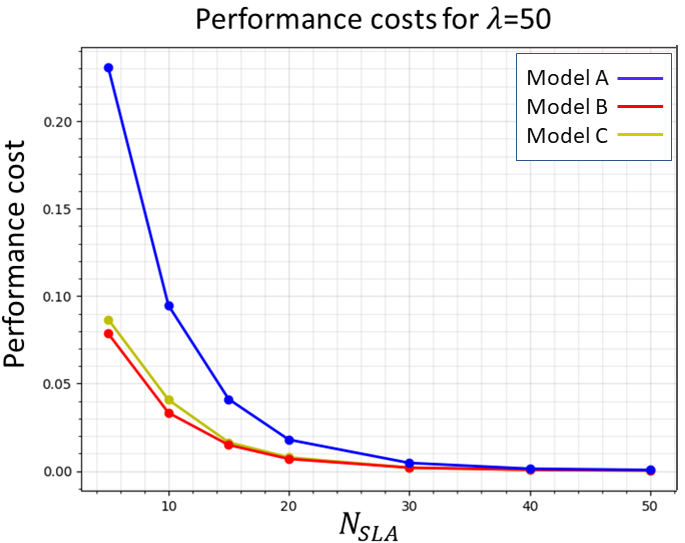}\par
\end{multicols}
\vspace{-0.6cm}
\begin{multicols}{2}
    \includegraphics[width=0.33\textwidth]{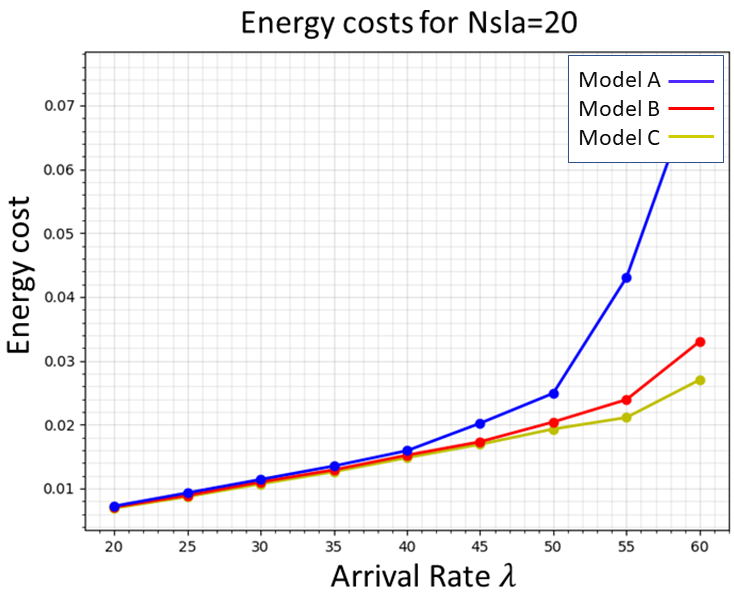}\par
    \includegraphics[width=0.33\textwidth]{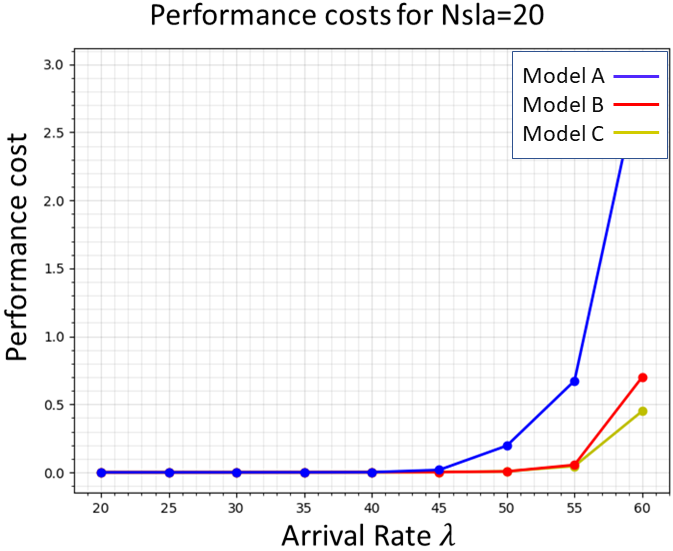}\par 
\end{multicols}
\caption{Energy and performance costs given fixed SLA or arrival rate}
\label{fig:eg_perf}
\end{figure}

\paragraph{Cost evaluation}
We display on Figure \ref{fig:eg_perf} an example ($K=8, B=60$) which exhibits the evolution of the financial cost per hour according to a given $N_{SLA}$, or to the load. We can see in the figure the impacts of such parameters regarding on performance and energy costs.

\vspace{-0.4cm}
\section{Conclusions and further work}
\label{sec:conclusion}
In this paper we have compared theoretically and numerically two different approaches to minimise the global cost integrating both performance and energy consumption in an auto scaling cloud system. 
The relevance of this study for a cloud provider is to provide an auto scaling 
policy very quickly to minimise its financial cost (around $38$ seconds for $128$ VMs and $800$ customers).
We exhibit that the best static heuristics are strongly outperformed by SMDP models and that hysteresis presents a significant improvement.
Our next steps will be to enlarge the model by including more general arrival process as burst traffic predicted from a measurement tool 
and second to build a cloud prototype to evaluate the performances of our policies in a real environment.  

\bibliographystyle{unsrt}  
\bibliography{biblioroadef}

\begin{thebibliography}{10}

\bibitem{report:oie}
{Observatoire de l’industrie Electrique}.
\newblock {\em Le cloud, les data centers et l’{é}nergie}, January 2017.

\bibitem{art:mastelic2015}
T.~{Mastelic} and I.~{Brandic}.
\newblock Recent trends in energy-efficient cloud computing.
\newblock {\em IEEE Cloud Computing}, 2:40--47, 2015.

\bibitem{art:lorido2014}
T.~Lorido-Botran, J.~Miguel-Alonso, and J.A. Lozano.
\newblock A review of auto-scaling techniques for elastic applications in cloud
  environments.
\newblock {\em J. Grid Computing}, 12:559–592, 2014.

\bibitem{art:benoit}
A.~Benoit, L.~Lef{\`e}vre, A.-C. Orgerie, and I.~Rais.
\newblock Reducing the energy consumption of large scale computing systems
  through combined shutdown policies with multiple constraints.
\newblock {\em Int. J. High Perform. Comput. Appl.}, 32(1):176--188, 2018.

\bibitem{art:asghari14}
N.M. Asghari, M.~Mandjes, and A.~Walid.
\newblock Energy-efficient scheduling in multi-core servers.
\newblock {\em Computer Networks}, 59:33--43, 2014.

\bibitem{art:serfozo}
M.Y. Kitaev and R.F. Serfozo.
\newblock {M}/{M}/1 queues with switching costs and hysteretic optimal control.
\newblock {\em Operations Research}, 47:310--312, 1999.

\bibitem{art:Lui99}
J.~C.~S. Lui and L.~Golubchik.
\newblock {Stochastic complement analysis of multi-server threshold queues with
  hysteresis}.
\newblock {\em Performance Evaluation}, 35:19--48, 1999.

\bibitem{conf:SPS15}
S.~Shorgin, A.~Pechinkin, K.~Samouylov, Y.~Gaidamaka, I.~Gudkova, and E.~Sopin.
\newblock Threshold-based queuing system for performance analysis of cloud
  computing system with dynamic scaling.
\newblock In {\em AIP Conference}, volume 1648, 2015.

\bibitem{art:Tou19}
T.~Tournaire, H.~Castel, E.~Hyon, and T.~Hoche.
\newblock Generating optimal thresholds in a hysteresis queue: a cloud
  application.
\newblock In {\em IEEE Mascots}, pages 283--294, 2019.

\bibitem{rap:aws}
Amazon.
\newblock {AWS} {A}uto {S}caling, 2018.

\bibitem{art:warsing2013}
D.P. Warsing, W.~Wangwatcharakul, and R.E. King.
\newblock Computing optimal base-stock levels for an inventory system with
  imperfect supply.
\newblock {\em Computers and Operations Research}, 40:2786--2800, 2013.

\bibitem{liv:song2013}
D.-P. Song.
\newblock {\em Optimal Control and Optimization of Stochastic Supply Chain
  Systems}.
\newblock Springer-Verlag, 2013.

\bibitem{art:orgerie}
M.~Kurpicz, A.-C. Orgerie, and A.~Sobe.
\newblock How much does a {VM} cost? energy-proportional accounting in
  {VM}-based environments.
\newblock In {\em PDP}, pages 651--658, 2016.

\bibitem{art:Krzy18}
J.~Krzywda, A.~Ali-Eldin, T.E. Carlson, P.-O. Ostberg, and E.~Elmroth.
\newblock Power-performance tradeoffs in data center servers: {DVFS}, {CPU}
  pinning, horizontal and vertical scaling.
\newblock {\em Future Generation Computer Systems}, pages 114--128, 2018.

\bibitem{art:adan2014}
I.J.B.F. Adan, V.G. Kulkarni, and A.C.C. van Wijk.
\newblock Optimal control of a server farm.
\newblock {\em INFOR: Information Systems and Operational Research},
  51(4):241--252, 2013.

\bibitem{Gandhi2010}
A.~Gandi, M.~Harchol-Balter, and I.~Adan.
\newblock Server farms with setup costs.
\newblock {\em Performance Evaluation}, 67(11):1123--1138, 2010.

\bibitem{art:mitrani2013managing}
I.~Mitrani.
\newblock Managing performance and power consumption in a server farm.
\newblock {\em Annals of Operations Research}, 202:121--134, 2013.

\bibitem{Artalejo2005}
J.~R. Artalejo, A.~Economou, and M.~J. Lopez-Herrero.
\newblock Analysis of a multiserver queue with setup times.
\newblock {\em Queueing Systems}, 51(1-2):53--76, 2005.

\bibitem{art:ardagna2014}
D.~Ardagna, G.~Casale, M.~Ciavotta, J.F. P{\'e}rez, and W.~Wang.
\newblock Quality-of-service in cloud computing: modeling techniques and their
  applications.
\newblock {\em JISA}, 5, 2014.

\bibitem{art:teghem85}
J.~Teghem.
\newblock Control of the service process in a queueing system.
\newblock {\em EJOR}, 23(2):141--158, 1986.

\bibitem{proc:YCNH11}
Z.~Yang, M.-H. Chen, Z.~Niu, and D.~Huang.
\newblock An optimal hysteretic control policy for energy saving in cloud
  computing.
\newblock In {\em GLOBECOM}, pages 1--5, 2011.

\bibitem{art:Lee2014}
N.~Lee and V.~G. Kulkarni.
\newblock Optimal arrival rate and service rate control of multi-server queues.
\newblock {\em Queueing Syst. Theory Appl.}, 76(1):37--50, 2014.

\bibitem{art:szarkowicz}
D.S. Szarkowicz and T.W. Knowles.
\newblock Optimal control of an {M}/{M}/{S} queueing system.
\newblock {\em Operations Research}, 33(3):644--660, 1985.

\bibitem{art:hipp84}
S.K. Hipp and U.D. Holzbaur.
\newblock Decision processes with monotone hysteretic policies.
\newblock {\em Operations Research}, 36(4):585--588, 1988.

\bibitem{art:NinoMora2019}
J.~Ni{ñ}o-Mora.
\newblock Resource allocation and routing in parallel multi-server queues with
  abandonments for cloud profit maximization.
\newblock {\em Computers and Operations Research}, 103:221--236, 2019.

\bibitem{art:randa2019}
A.C. Randa, M.K. Dogru, C.~Iyigun, and U.~Özen.
\newblock Heuristic methods for the capacitated stochastic lot-sizing problem
  under the static-dynamic uncertainty strategy.
\newblock {\em Computers and Operations Research}, 109:89--101, 2019.

\bibitem{art:kranenburg}
A.A. Kranenburg and G.J. van Houtum.
\newblock Cost optimization in the ({S}-1, {S}) lost sales inventory model with
  multiple demand classes.
\newblock {\em Oper. Res. Lett.}, 35(4):493--502, 2007.

\bibitem{art:braglia2019}
M.~Braglia, D.~Castellano, L.~Marrazzini, and D.~Song.
\newblock A continuous review, ({Q}, r) inventory model for a deteriorating
  item with random demand and positive lead time.
\newblock {\em Computers and Operations Research}, 109:102--121, 2019.

\bibitem{art:Wu2014}
C.-H. Wu, W.-C. Lee, J.-C. Ke, and T-H. Liu.
\newblock Optimization analysis of an unreliable multi-server queue with a
  controllable repair policy.
\newblock {\em Computers and Operations Research}, 49:83–96, 2014.

\bibitem{art:ibekeilson95}
O.C. Ibe and J.~Keilson.
\newblock Multi-server threshold queues with hysteresis.
\newblock {\em Performance Evaluation}, 21:185--213, 1995.

\bibitem{conf:epewKandi17}
M.M. Kandi, F.~A{\"{\i}}t{-}Salaht, H.~Castel{-}Taleb, and E.~Hyon.
\newblock Analysis of performance and energy consumption in the cloud.
\newblock In {\em {EPEW}}, pages 199--213, 2017.

\bibitem{liv:put1994}
M.L. Puterman.
\newblock {\em Markov Decision Processes: Discrete Stochastic Dynamic
  Programming}.
\newblock Wiley, 1994.

\bibitem{conf:ajmMarmote}
A.~Jean-Marie.
\newblock marmote{C}ore: a {M}arkov modeling platform.
\newblock In {\em VALUETOOLS}, pages 60--65, 2017.

\bibitem{rap:aws2}
Amazon.
\newblock Amazon {EC2} pricing, 2019.

\bibitem{art:orgerieWL}
D.~Guyon, A-C. Orgerie, C.~Morin, and D.~Agarwal.
\newblock Involving users in energy conservation: A case study in scientific
  clouds.
\newblock {\em International Journal of Grid and Utility Computing}, 2018.

\end{thebibliography}

\appendix

\section{Resources management examples}
\label{app:example}
We detail now the management resources.
\begin{example}
We display an example for all cases according to Definition \ref{def:MDPscenarios}, to illustrate the difference of optimal solutions between both approaches. Solutions returned by MDP algorithms are displayed without the shift operation.
We take the following costs values $C_a = C_d = 2, ~ C_h = C_s = 5, ~ C_r = 10$ with $K=16, ~ B=100$ scenario.

\paragraph{Medium arrival case}
Consider the following queueing parameters: $\lambda = 500$, $\mu=100$.

The MDP algorithm provides the following solution:
$C^* = 25,239 \mbox{; and}$ 

$L^* = [2,4,5,7,8,10,11,13,15,16,18,19,21,23,24];$

$l^* = [0,1,3,4,5,6,7,8,9,10,11,12,13,14,15].$

The heuristic in the (MC) model provides the following solution:
$C^* = 25,239 \mbox{; and :}$ 

$F^* = [1,3,4,6,7,9,10,12,14,15,17,18,20,22,23];$

$R^* = [1,2,4,5,6,7,8,9,10,11,12,13,14,15,16].$

We can observe here that the two solutions are identical since (MDP) algorithm finds thresholds in all levels.

\paragraph{Low arrival case}
Consider the following queueing parameters: $\lambda = 50$, $\mu=100$.

The MDP algorithm provides the following solution:
$C^* = 9,99323 \mbox{; and :}$ 

$L^* = [12,26,39,52,64,77,90,\infty,\infty,\infty,\infty,\infty,\infty,\infty,\infty];$

$l^* = [0,1,2,3,4,5,6,7,8,9,10,11,12,13,14].$

The heuristic in the (MC) model provides the following solution:
$C^* = 9,99375 \mbox{; and :}$ 

$F^* = [11,25,28,39,51,58,70,72,88,94,95,96,97,98,99];$

$R^* = [1,2,3,4,5,6,7,8,9,10,11,12,13,14,15].$

Since the load in the system is very low, the (MDP) algorithm does not need to turn ON all virtual resources. However, we observe that (MC) approach still activate all levels but tries to activate as late as possible to compensate that it does not need to turn ON more virtual resources. It induces a small supplementary costs since these activations will 'never' happen in this system.

\paragraph{High arrival case}
Consider the following queueing parameters: $\lambda = 1000$, $\mu=100$.

The MDP algorithm provides the following solution:
$C^* = 119,756 \mbox{; and :}$ 

$L^* = [2,3,4,5,6,7,9,10,11,13,16,24,56,69,82];$

$l^* = [0,0,0,\ldots,0,0,0,10,12,14].$

The heuristic in the (MC) model provides the following solution:
$C^* = 220,879 \mbox{; and :}$ 

$F^* = [2,4,6,8,11,13,16,21,27,37,48,61,75,88,96];$

$R^* = [1,2,3,4,5,6,7,8,11,16,22,26,29,30,34].$

We can observe here that (MC) approach still define deactivation thresholds in any levels (due to its model and isotone hysteresis constraints) leading to a higher mean cost.
\end{example}

\section{Complements on SMDP anaysis}
\label{app:analysis}.
We give now some examples to illustrate the theoretical analysis of the SMDP.

\subsection{Multichain MDP example}

\begin{example}
We show a scenario where the MDP policy $q$ will generate a multichain Markov chain where we can observe several recurrent classes. This example is represented on Figure \ref{fig:multichain}.

\begin{figure}[hbtp]
\centering
    \includegraphics[scale=0.5]{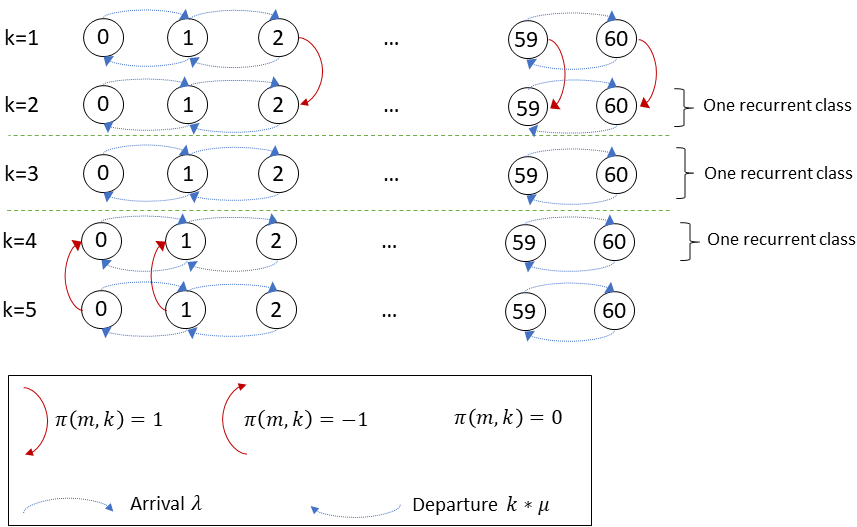}
    \caption{Example of an isotone hysteresis multichain MDP}
    \label{fig:multichain}
\end{figure}
\end{example}

\subsection{Temporal behaviour of two approaches}
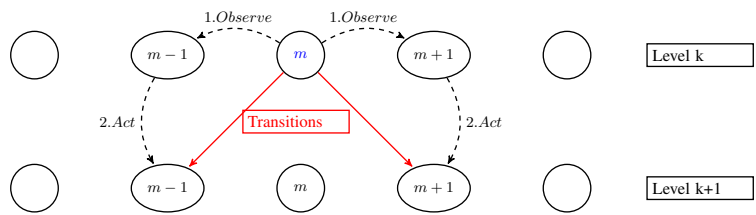
\begin{figure}[hbtp]
\centering
\caption{Behaviour of the agent in the Markov chain Model}
\label{MC_transitions}
\resizebox{10cm}{!}{
\begin{tikzpicture}[->,>=stealth',shorten >=1pt,auto,node distance=2.8cm,
  thick,main node/.style={ellipse,fill=white!10,draw,minimum size=1cm,font=\sffamily\small\bfseries}]
  \node[main node] (E1) {};
  \node[main node] (E2) [right of=E1] {$m-1$};
  \node[main node] (E3) [right of=E2] {$\textcolor{blue}{m}$};
  \node[main node] (E4) [right of=E3] {$m+1$};
  \node[main node] (E5) [right of=E4] {};
 
  \node[main node] (E6) [below of=E1] {};
  \node[main node] (E7) [right of=E6] {$m-1$};
  \node[main node] (E8) [right of=E7] {$m$};
  \node[main node] (E9) [right of=E8] {$m+1$};
  \node[main node] (E10) [right of=E9] {};
  
  \path[every node/.style={font=\sffamily\small}]
  
   (E3) edge [bend left=330,above,dashed] node {$1. Observe$} (E2)
   (E2) edge [bend left=330,left,dashed] node {$2. Act$} (E7)
   
   (E3) edge [bend right=330,dashed] node {$1. Observe$} (E4)
   (E4) edge [bend right=330,dashed] node {$2. Act$} (E9)
   
   (E3) edge [color=red] node {} (E7)
   (E3) edge [color=red] node {} (E9);
   
   \node[draw,text width=2cm,color=red] at (5.5,-1.4) {\textcolor{red}{Transitions}};
   \node[draw,text width=2cm] at (14,0) {Level k};
   \node[draw,text width=2cm] at (14,-2.8) {Level k+1};
   
     \end{tikzpicture}
     }
\end{figure}

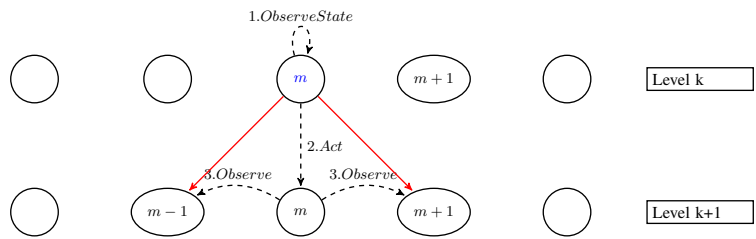
\begin{figure}[hbtp]
\centering
\caption{Behaviour of the agent in the MDP Model}
\label{MDP_transitions}
\resizebox{10cm}{!}{
\begin{tikzpicture}[->,>=stealth',shorten >=1pt,auto,node distance=2.8cm,
  thick,main node/.style={ellipse,fill=white!10,draw,minimum size=1cm,font=\sffamily\small\bfseries}]
  
  \node[main node] (E1) {};
  \node[main node] (E2) [right of=E1] {};
  \node[main node] (E3) [right of=E2] {$\textcolor{blue}{m}$};
  \node[main node] (E4) [right of=E3] {$m+1$};
  \node[main node] (E5) [right of=E4] {};
 
  \node[main node] (E6) [below of=E1] {};
  \node[main node] (E7) [right of=E6] {$m-1$};
  \node[main node] (E8) [right of=E7] {$m$};
  \node[main node] (E9) [right of=E8] {$m+1$};
  \node[main node] (E10) [right of=E9] {};
  
  \path[every node/.style={font=\sffamily\small}]
  
   (E3) edge [loop above,dashed] node {$1. Observe State$}
   (E3) edge [dashed] node {$2. Act$} (E8)
   (E8) edge [bend right=330,dashed] node {$3. Observe$} (E9)
   (E8) edge [bend left=330, dashed,above] node {$3. Observe$} (E7)
   
   (E3) edge [color=red] node {} (E7)
   (E3) edge [color=red] node {} (E9);

   \node[draw,text width=2cm] at (14,0) {Level k};
   \node[draw,text width=2cm] at (14,-2.8) {Level k+1};
   
     \end{tikzpicture}
     }
\end{figure}

We display two figures showing the different temporal behaviour between the two approaches: Markov chain and MDP.
The first figure (Figure \ref{MC_transitions}) explains the temporal transition in the Markov chain approach while the second figure (Figure \ref{MDP_transitions}) explains the temporal transition for MDP.

\subsection{Transformation of MDP policy into a threshold policy}
We display an example of generated Markov chain by a policy $q$ of the MDP agent. We show in this figure how we can transform the MDP policy into a threshold policy. 
\begin{example}
\label{exemple:transformThreshold}
Example of the transformation of a MDP policy into a threshold policy with 3 virtual machines. This case is illustrated on Figure \ref{fig:transformPol}.
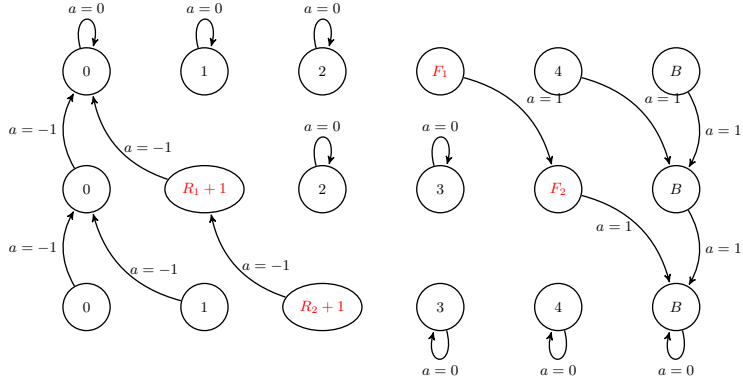
\begin{figure}[hbtp]
\centering
\resizebox{10cm}{!}{
\begin{tikzpicture}[->,>=stealth',shorten >=0.8pt,auto,node distance=2.5cm,
  thick,main node/.style={ellipse,fill=white!10,draw,minimum size=1cm,font=\sffamily\small\bfseries}]
  
  \node[main node] (E1) {$0$};
  \node[main node] (E2) [right of=E1] {$1$};
  \node[main node] (E3) [right of=E2] {$2$};
  \node[main node] (E4) [right of=E3] {\textcolor{red}{$F_1$}};
  \node[main node] (E5) [right of=E4] {$4$};
  \node[main node] (E6) [right of=E5] {$B$};
 
  \node[main node] (E7) [below of=E1] {$0$};
  \node[main node] (E8) [right of=E7] {\textcolor{red}{$R_1 + 1$}};
  \node[main node] (E9) [right of=E8] {$2$};
  \node[main node] (E10) [right of=E9] {$3$};
  \node[main node] (E11) [right of=E10] {\textcolor{red}{$F_2$}};
  \node[main node] (E12) [right of=E11] {$B$};
  
  \node[main node] (E13) [below of=E7] {$0$};
  \node[main node] (E14) [right of=E13] {$1$};
  \node[main node] (E15) [right of=E14] {\textcolor{red}{$R_2 + 1$}};
  \node[main node] (E16) [right of=E15] {$3$};
  \node[main node] (E17) [right of=E16] {$4$};
  \node[main node] (E18) [right of=E17] {$B$};
  
  \path[every node/.style={font=\sffamily\small}]
  
   (E1) edge [loop above] node {$a=0$} (E1)
   (E2) edge [loop above] node {$a=0$} (E2)
   (E3) edge [loop above] node {$a=0$} (E3)
   (E4) edge [bend right=330] node {$a=1$} (E11)
   (E5) edge [bend right=330] node {$a=1$} (E12)
   (E6) edge [bend right=330] node {$a=1$} (E12)
   
   (E7) edge [bend left] node {$a=-1$} (E1)
   (E8) edge [bend left, right] node {$a=-1$} (E1)
   (E9) edge [loop above] node {$a=0$} (E9)
   (E10) edge [loop above] node {$a=0$} (E10)
   (E11) edge [bend right=330, left] node {$a=1$} (E18)   
   (E12) edge [bend right=330] node {$a=1$} (E18)
   
   (E13) edge [bend left] node {$a=-1$} (E7)
   (E14) edge [bend left,right] node {$a=-1$} (E7)
   (E15) edge [bend left,right] node {$a=-1$} (E8)
   (E16) edge [loop below] node {$a=0$} (E9)
   (E17) edge [loop below] node {$a=0$} (E16)   
   (E18) edge [loop below] node {$a=0$} (E17);
   		
     \end{tikzpicture}
     }
     \caption{Example of the threshold policy transformation in a low scale scenario}
     \label{fig:transformPol}
\end{figure}

The optimal threshold vector in this example would be: $L^* = [3,4] \mbox{ and } l^* = [0,1]$
\end{example}

\section{Experimental results for concrete scenarios}
\label{app:experiments}
\begin{figure}[hbtp]
\begin{multicols}{2}
    \includegraphics[scale=0.4]{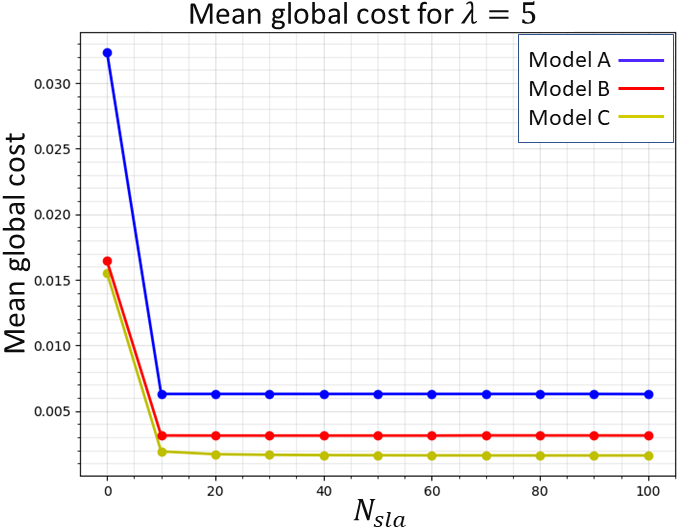}\par
    \includegraphics[scale=0.4]{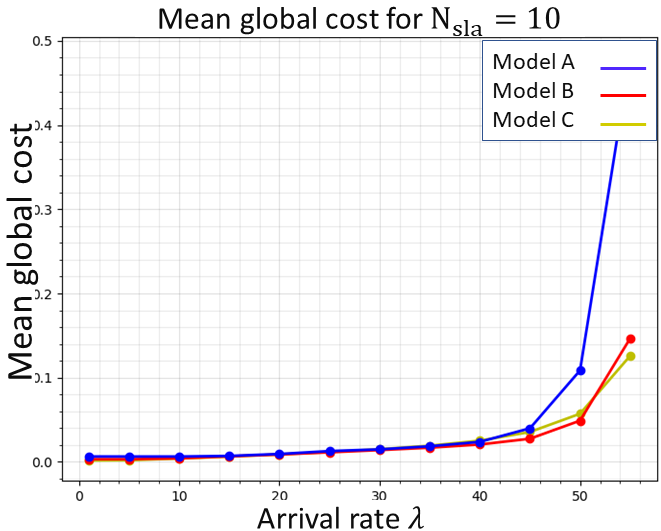}\par 
    \end{multicols}
\begin{multicols}{2}
    \includegraphics[scale=0.4]{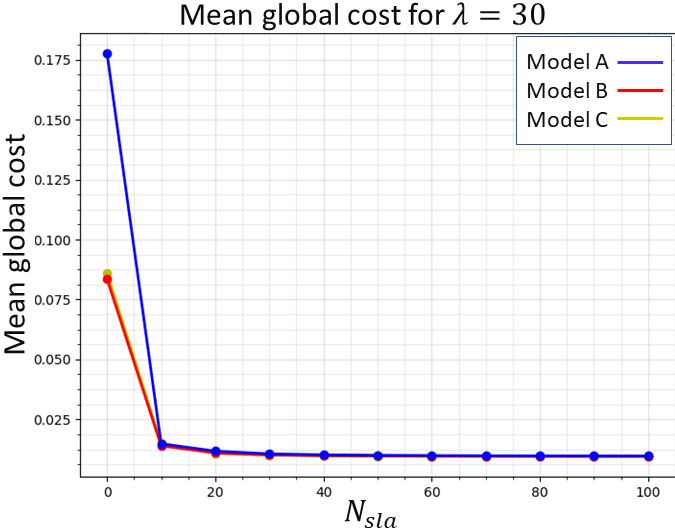}\par
    \includegraphics[scale=0.4]{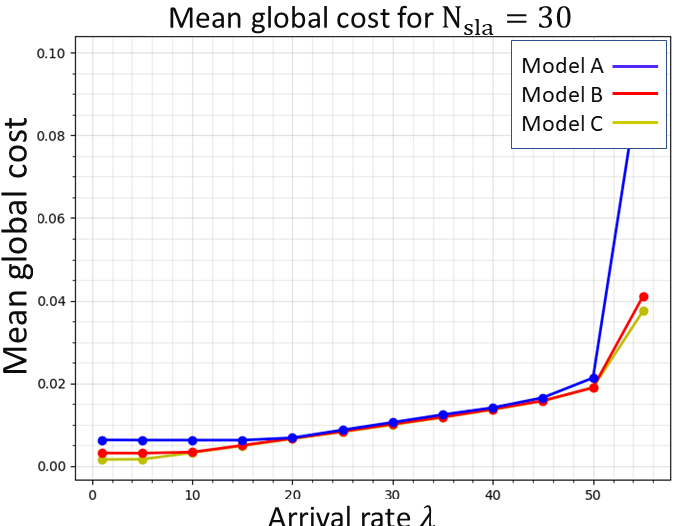}\par
\end{multicols}
\begin{multicols}{2}
    \includegraphics[scale=0.4]{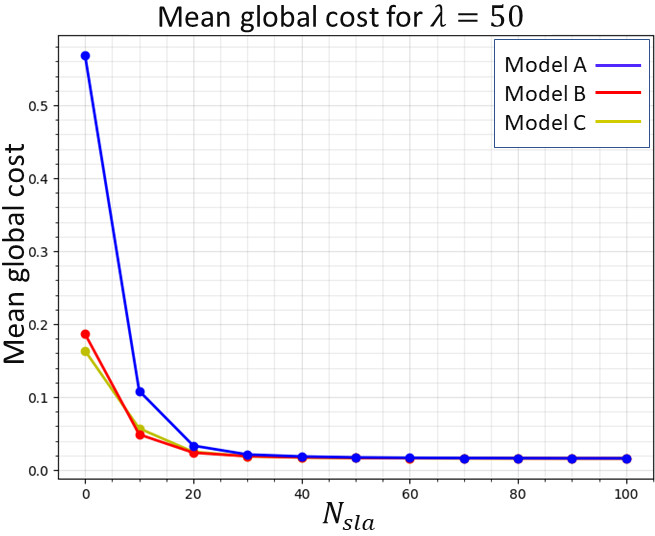}\par
    \includegraphics[scale=0.4]{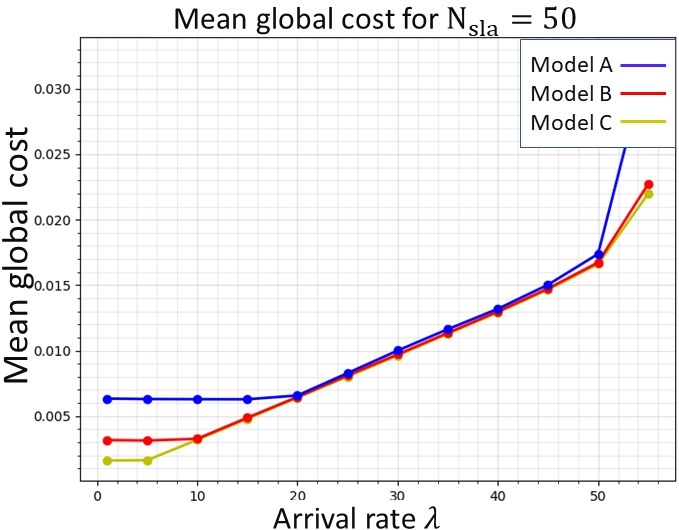}\par
\end{multicols}
\caption{Mean global cost for the three models given a fixed Nsla or a fixed lambda}
\label{fig:meancostconcrete}
\end{figure}

\begin{figure}[hbtp]
\begin{multicols}{2}
    \includegraphics[scale=0.4]{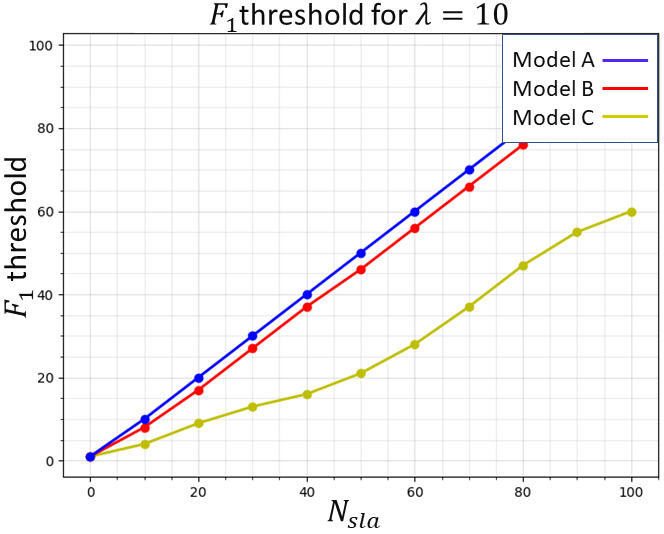}\par 
    \includegraphics[scale=0.4]{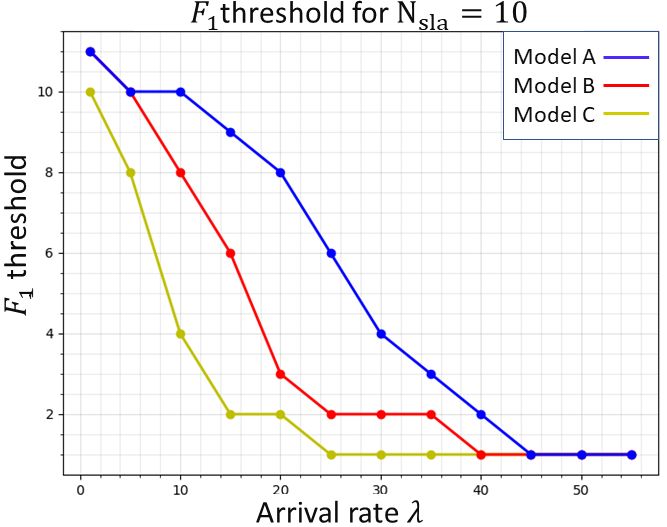}\par 
    \end{multicols}
    
\begin{multicols}{2}

    \includegraphics[scale=0.4]{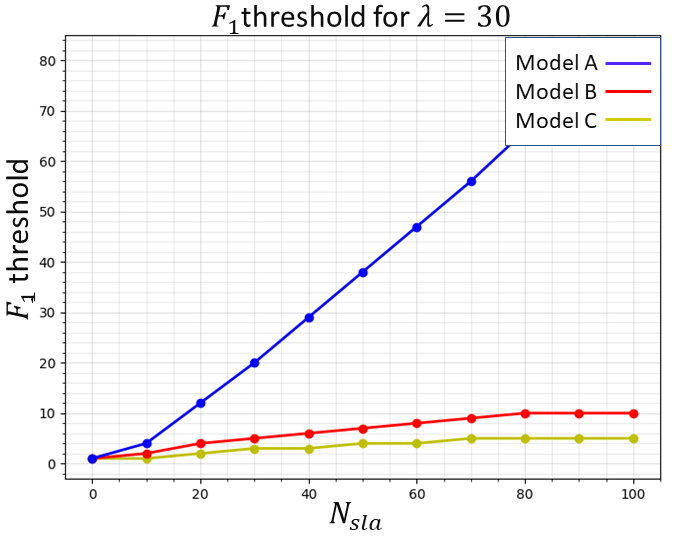}\par
    \includegraphics[scale=0.4]{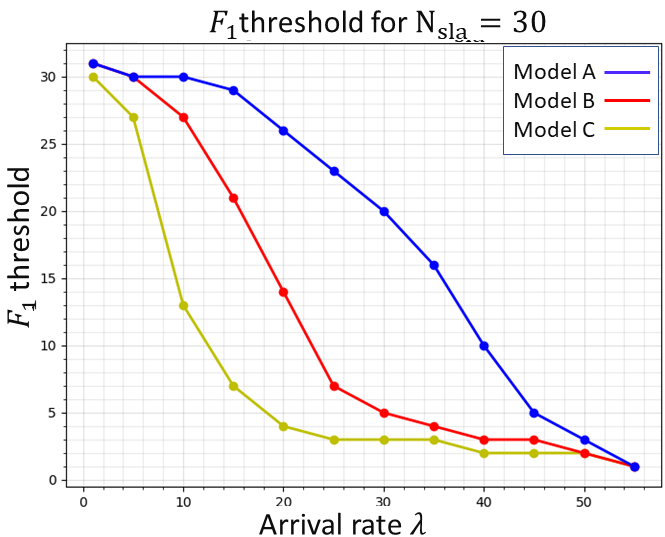}\par
\end{multicols}

\begin{multicols}{2}
    \includegraphics[scale=0.4]{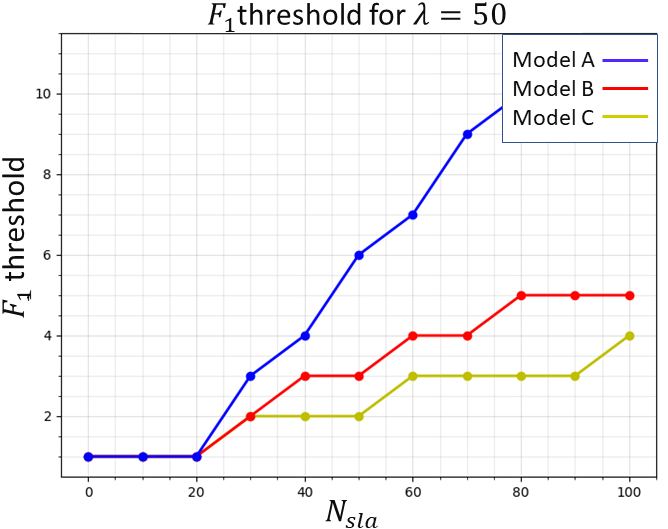}\par
    \includegraphics[scale=0.4]{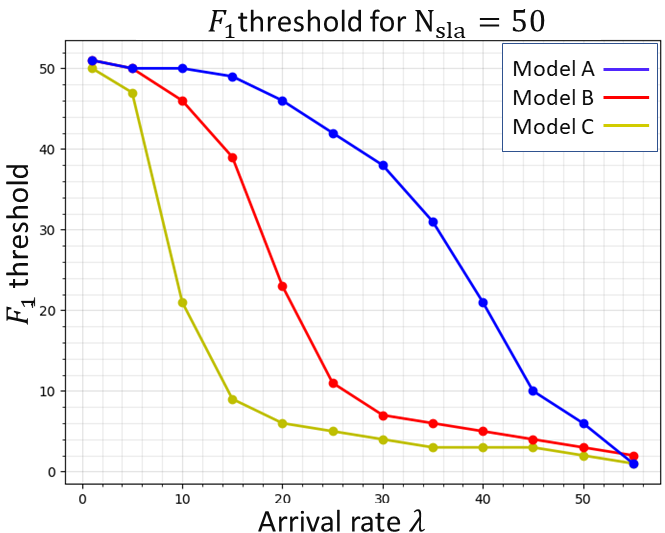}\par
\end{multicols}

\caption{F1 thresholds for the three models given a fixed Nsla or a fix lambda}
\label{fig:F1concrete}
\end{figure}

We display several supplementary experimental (Section 7.3) results 
for concrete scenarios by comparing the behaviour of first activation thresholds and global costs when we take different parameterisation of the system (arrival rate $\lambda$, SLA metric $N_{sla}$).
These are given in the figures later. Figure \ref{fig:meancostconcrete} represents the mean global cost while Figure \ref{fig:F1concrete} gives the evolution of the first threshold.

\section{Pseudo-Codes of algorithms}
We present the pseudo-codes of heuristics described in section 4.3.
To describe the different algorithms, we define the following notations: 

- \textit{solve$[F,R]$}: function that, for a given vector of thresholds, generates its associated Markov chain, computes its stationary distribution with the aggregation method then computes the mean cost $\overline{C}_{[F,R]}^{\pi}$ .

- $C^*$: the solution of the optimisation problem and  $[F,R]^*$ the optimal vector of thresholds. 

\begin{algorithm}[hbtp]
\DontPrintSemicolon
 
  \KwInput{$\lambda, \mu, B, K, C_a, C_d, C_h, C_s, C_r$,  the cloud system parameters}
  \KwOutput{$[F,R]^*, C^*$}
  \textbf{Initialise} a vector of thresholds $[F,R]_0$
  
  Apply \textbf{solve}$[F,R]$ for the current vector

  \While{improvement=TRUE}
    {   \tcc{Activation Thresholds}
        
        \For{$k \in [1,K-1]$}
            { 
                \For{$F_k$ respecting constraints}
                    {   
                        \textbf{solve}$[F,R]$ and \textbf{store} the vector if improvement
                    }
            }
            
        \hrulefill
        
        \tcc{Deactivation Thresholds}
        
        \For{$k \in [1,K-1]$}
            { 
                \For{$R_k$ respecting constraints}
                    {   
                        \textbf{solve}$[F,R]$ and \textbf{store} the vector if improvement
                    }
            }
    }
\label{algo:bpl}
\caption{Best per level}
\end{algorithm}

\begin{algorithm}[hbtp]
\DontPrintSemicolon
 
  \KwInput{$\lambda, \mu, B, K, C_a, C_d, C_h, C_s, C_r$,  the cloud system parameters}
  \KwOutput{$[F,R]^*, C^*$}
  \textbf{Initialise} a vector of thresholds $[F,R]_0$
  
  Apply \textbf{solve}$[F,R]$ for the current vector
  
  \While{improvement=TRUE}
  { \tcc{We generate $\mathbf{V}$ by giving +1 or -1 to each elements of the vector, $\text{card}(\mathbf{V}) = 4(K-1)$} 
    \textbf{Generate} the neighbourhood $\mathbf{V}$ of $[F,R]$ 
    
    \For{$[F,R] \in \mathbf{V}$}
        {
            \textbf{solve}$[F,R]$ and \textbf{store} the vector if improvement
        }
  }
\label{algo:nls}
\caption{Neighbourhood local search}
\end{algorithm}

\begin{algorithm}[hbtp]
\DontPrintSemicolon
 
  \KwInput{$\lambda, \mu, B, K, C_a, C_d, C_h, C_s, C_r$,  the cloud system parameters}
  \KwOutput{$[F,R]^*, C^*$}
  \textbf{Initialise} a vector of thresholds $[F,R]_0$
  
  Apply \textbf{solve}$[F,R]$ for the current vector

  \While{improvement=TRUE}
    {  \tcc{Activation Threshold 1st Level}
        
        First level k=1 
        
        \For{$m \in [1,B]$}
        {
            Compute $\phi^A_1(m)$
            
            When $\phi^A_1(m)$ > 0 ,BREAK, then $F_1$ = m if constraints respected
        }

      \tcc{Activation thresholds level k}
        
      \For{$k \in [2,K-1]$}
        {
            \For{$ m \in [F_{k-1}+1,B]$}
            {
              Compute $\phi^A_k(m)$
              
              When $\phi^A_k(m)$ > 0 ,BREAK, then $F_k$ = m if constraints respected
            }
        }

        \tcc{Deactivation Threshold 1st Level}
        
        First level k=1 
        
        \For{$m \in [0,F_k - 1]$}
        {
            Compute $\phi^D_1(m)$
            
            When $\phi^D_1(m) \geq 0$ ,BREAK, then $R_1$ = m if constraints respected
        }

      \tcc{Deactivation thresholds level k}
        
      \For{$k \in [2,K]$}
        {
            \For{$m \in [R_{k-1}+1,F_k - 1]$}
            {
              Compute $\phi^D_k(m)$
              
              When $\phi^D_k(m) \geq 0$ ,BREAK, then $R_k$ = m if constraints respected
            }
        }
    }
\label{algo:mmk}
\caption{M/M/K/B approximation}
\end{algorithm}

\end{document}